\documentstyle{article}     

\newtheorem{th}{Theorem}
\newtheorem{ax}{Axiom}
\newtheorem{lm}{Lemma}
\newtheorem{df}{Definition}
\newtheorem{pr}{Proposition}
\newtheorem{cl}{Corollary}
\newtheorem{as}{Assumption}
\newtheorem{qu}{Question}
\newcommand{\bth}{\begin{th}\hspace{-5pt}{\bf .} \ }
\newcommand{\eth}{\end{th}}
\newcommand{\bax}{\begin{ax}\hspace{-5pt}{\bf .} \ }
\newcommand{\eax}{\end{ax}}
\newcommand{\blm}{\begin{lm}\hspace{-5pt}{\bf .} \ }
\newcommand{\elm}{\end{lm}}
\newcommand{\bdf}{\begin{df}\hspace{-5pt}{\bf .} \ }
\newcommand{\edf}{\end{df}}
\newcommand{\bpr}{\begin{pr}\hspace{-5pt}{\bf .} \ }
\newcommand{\epr}{\end{pr}}
\newcommand{\bcl}{\begin{cl}\hspace{-5pt}{\bf .} \ }
\newcommand{\ecl}{\end{cl}}
\newcommand{\bas}{\begin{as}\hspace{-5pt}{\bf .} \ }
\newcommand{\eas}{\end{as}}
\newcommand{\bqu}{\begin{qu}\hspace{-5pt}{\bf .} \ }
\newcommand{\equ}{\end{qu}}
\newcommand{\bpf}{\noindent {\it Proof}\hspace{0.1truecm}: }
\newcommand{\epf}{\hfill${\Box}$\par\vspace{2.5mm}\noindent}
\newcommand{\bit}{\begin{itemize}}
\newcommand{\eit}{\end{itemize}\par\noindent}
\newcommand{\beq}{\begin{equation}}
\newcommand{\eeq}{\end{equation}\par\noindent}
\newcommand{\beqa}{\begin{eqnarray*}}
\newcommand{\eeqa}{\end{eqnarray*}\par\noindent}
\newcommand{\beqn}{\begin{eqnarray}}
\newcommand{\eeqn}{\end{eqnarray}\par\noindent}

\newcommand{\ela}{\end{array}\right.}
\newcommand{\bra}{\left\{\begin{array}{r}}
\newcommand{\era}{\end{array}\right.}

\newcommand{\ben}{\begin{enumerate}}
\newcommand{\een}{\end{enumerate}\par\noindent}

\begin{document} 
\
\par\medskip\par\noindent
\centerline{\large{\bf QUANTUM LOGIC IN INTUITIONISTIC PERSPECTIVE}}
\par\medskip\par\noindent
\centerline{\normalsize{B{\scriptsize OB} C{\scriptsize OECKE}}}
\par\medskip\par\noindent
\centerline{\small{Free University of Brussels, Department of Mathematics,}}
\par\noindent
\centerline{\small{Pleinlaan 2, B-1050 Brussels\,; bocoecke@vub.ac.be\,;}}
\par\smallskip\par\noindent
\centerline{\small{\sl and}}
\par\smallskip\par\noindent
\centerline{\small{Imperial College of Science, Technology \&
Medicine, Theoretical Physics Group,}} 
\par\noindent 
\centerline{\small{The Blackett Laboratory, South Kensington, London SW7 2BZ\,.}}
\par\bigskip\par\noindent 
\centerline{\small{{\it Current address:}}}

\smallskip\noindent
\centerline{\small{University of Oxford, Computing 
Laboratory\,,  Wolfson Building, Parks Road,}} 

\noindent 
\centerline{\small{Oxford, OX1 3QD, UK\,; e-mail:
coecke@comlab.ox.ac.uk\,.}} 
\begin{abstract}
\noindent
In their seminal paper Birkhoff and von Neumann revealed the
following dilemma\,: ``... whereas for
logicians the orthocomplementation properties of negation were the
ones least able to withstand a critical
analysis, the study of mechanics points to the distributive
identities as the weakest link in the algebra
of logic.''~~In this paper we eliminate this dilemma, providing a way
for maintaining both.
Via the introduction of the ``missing'' disjunctions in the lattice of
properties of a physical system while inheriting the meet as a
conjunction we obtain a complete Heyting
algebra of propositions on physical properties. In particular there
is a bijective correspondence
between property lattices and propositional lattices equipped with a
so called {\it operational
resolution}\,, an operation that exposes the properties on the level
of the propositions. If the
property lattice goes equipped with an orthocomplementation, then
this bijective correspondence can be
refined to one with propositional lattices equipped with an {\it
operational complementation}\,, as such establishing the claim made above.
Formally
one rediscovers via physical and logical considerations as such
respectively a specification and a
refinement of the purely mathematical result by Bruns and Lakser
(1970) on injective hulls of
meet-semilattices. From our representation we can derive  a truly
intuitionistic functional implication
on property lattices, as such confronting claims made in previous
writings on the matter. We also make a
detailed analysis of disjunctivity vs.~distributivity and finitary
vs.~infinitary conjunctivity, we
briefly review the Bruns-Lakser construction and indicate some
questions which are left open.
\end{abstract}

\bigskip\noindent
{\bf 0. CONTENT}

{\footnotesize
\medskip\noindent
1. INTRODUCTION 

\smallskip\noindent
2. FORMAL AND METHODOLOGICAL TOOLS

\smallskip
2.1. Logical significance of property lattices

\smallskip
2.2. Bruns-Lakser distributive hulls for complete lattices

\smallskip\noindent
3. MAIN RESULTS

\smallskip
3.1. Complete lattices and operational resolution

\smallskip
3.2. Complete ortholattices and operational complementation

\smallskip\noindent
4. FURTHER ANALYSIS AND OPEN PROBLEMS

\smallskip
4.1. Characterization of disjunctivity

\smallskip
4.2. Finitely conjunctive infima

\smallskip\noindent
5. SUMMARY, CONCLUSION AND PERSPECTIVES

\smallskip\noindent
APPENDIX: IMPLICATION VIA ADJUNCTION 
}

\bigskip\noindent
{\bf 1. INTRODUCTION} 

\medskip\noindent
In their seminal paper Birkhoff and von
Neumann (1936) observe that the lattice of closed subspaces of a
Hilbert space retains a
number of the familiar features of Boolean algebras (which constitute
the semantics of
classical propositional logic), namely, it is orthocomplemented and
hence satisfies the De
Morgan laws. However, the distributive law fails. Confronting the then ongoing
tendencies towards intuitionistic logic [Birkhoff and von Neumann
(1936) p.839]\,:
\begin{quote}
``The models for propositional calculi [of physically significant
statements in quantum
mechanics] are also interesting from the standpoint of pure logic.
Their nature is determined by quasi-physical and technical reasoning,
different from the
introspective and philosophical considerations which have to guide
logicians hitherto [\,...\,]
whereas logicians  have usually assumed that [the
orthocomplementation] properties
L71-L73 of negation were the ones least able to withstand a critical
analysis, the study of
mechanics points to the distributive identities L6 as the weakest
link in the algebra
of logic.''
\end{quote}
they point at a fundamental
difference between Heyting algebras (the semantics of
intuitionistic propositional logic) and orthomodular lattices (the
``usual'' semantics of
quantum logic) when  viewed as generalizations of Boolean algebra.
This seems to enforce a
{\it dilemma}\, with respect to logical considerations on
propositions attributed to physical
systems.  It is probably fair to say that due to this dilemma,
quantum logic became a
strictly separated domain of mathematics that had no essential impact
on traditional
fields of logic. Moreover, most attempts to provide a logical syntax
for discussing
physical properties, e.g., Hardegree (1979) and Kalmbach (1993), knew
serious criticism
(and definitely not always unjust), e.g., the arguments in Goldblatt
(1984), Malinowski (1990)
and Moore (1993).  In particular we do want to point in this context
to the failure to equip
quantum logic with a satisfactory internal implication operation.

\smallskip
However, we will show in this paper that both motivations, i.e., the
physical one encoding a
non-distributive orthocomplemented lattice and the logical
intuitionistic one encoding
a distributive pseudocomplemented lattice, are not incompatible but
motivate a distinction
between the {\sl physical properties} themselves and {\sl logical
propositions on physical
properties}.  In particular we will encode all aspects within one
mathematical object, namely
a complete Heyting algebra equipped with an additional operation, the
{\sl operational
resolution}. We will motivate these claims and constructions using
the {\sl operational} {\it
methodology}\,\footnote{By some people considered as a {\it
doctrine}, including one of the
fathers, namely Piron himself\,; dixit Piron\,: ``Les Coeckeries et
les Moorismes ne sont pas
des Pironeries''.  In our view, the operational methodology allows to
communicate and refine certain
insights, independent on the reader's personal view on physics, and
contributes either in providing an image, an understanding or a model,
this depending on the reader's personal taste.}
for quantum logic, which was already implicitly indicated in Birkhoff
and von Neumann (1936), but got
only truly established in Jauch and Piron (1969) and further
developed and refined in Piron (1976), Aerts
(1982), Moore (1999) and Coecke, Moore and Smets (2001a,b). This
methodology relates properties of a
physical system to definite experimental projects in part to provide
an answer to [Birkhoff and von
Neumann (1936) p.839]\,:
\begin{quote}
``What experimental meaning can one attach to the meet and join of two given
experimental propositions [on quantum systems]?''
\end{quote}
but also to motivate a common framework to discuss both classical and
quantum systems,
and understand their ontological and epistemological differences.
However, we feel that even
if one does not fully subscribe to this methodology, most, and in
particular all essential aspects of this
paper still hold.  For example, one ingredient of the methodology
consists of {\it  proving}\, that the
lattice of properties of a physical system should be taken complete,
i.e., any subset of it has
a greatest lower bound and a least upper bound. However, if one
chooses to think of a
property lattice as having only finite meets, or if one rather has
analytical or
probabilistic inspirations, as such {\it preferring}\, it to be
$\sigma$-complete, the construction and corresponding interpretation
of greatest lower bounds
for arbitrary (large) sets of properties provides a way to think about
completion of this property lattice, and how one should manipulate
this extension.
In particular is Section 4.2 of this paper devoted to property lattices in which
only finite meets are considered as conjunctions.

\smallskip
A striking
fact of the mathematics applied in this paper is indeed that the
assumption on preservation of
finite meets in the considered representation automatically ensures
preservation of all
infinitary meets as well.
About this mathematics, the in this paper proposed representation for the
properties of a physical system within a complete Heyting algebra (of
logical propositions on
these properties) equipped with a particular kind of closure operator (the
operational resolution), and
which will be motivated by logical reflection on primitive
operational physical notions,
this representation actually mimics a purely mathematically motivated
result of Bruns and
Lakser (1970), also independently found by Horn and Kimura (1971),
namely proving the
existence and characterizing the injective hulls in the
category of meet-semilattices.
We will specify this result for complete lattices (and show
that this specification
works), thereby moulding it towards our particular needs.

\smallskip
Concluding this introduction, in order to
substantiate our claim at the beginning of the previous paragraph we
will specify this
representation for complete ortholattices, as such revealing the
physical notion of an {\sl
operational complementation}, a pseudo-orthocomplementation that has
the operational
resolution as its square\,: This operational complementation will
then be the operation that
recaptures the orthocomplementation of the properties as an
additional operation on the
complete Heyting algebra of the logical propositions on these properties.
Since in this representation it is the collection of propositions
that goes equipped with an
{\sl internal intuitionistic implication operation}, ``implication
for physical properties''
should be envisioned as an external operation that assigns
propositions to pairs of properties.

\vfill\eject 
\noindent
{\bf 2. FORMAL AND METHODOLOGICAL TOOLS}

\medskip\noindent
First we provide and discuss the required formal and methodological tools.

\bigskip
\noindent
{\bf 2.1. Logical significance of property lattices}

\medskip\noindent
Let us briefly survey the fragment of the above mentioned operational
methodology that we will
employ in this paper\,; we refer to Moore (1999) and Coecke, Moore
and Smets (2001a,b) for the
most recent overview.\footnote{In particular contain the latter two
of these papers a
critical analysis of the methodology as part of situating it within
philosophy of
science as a whole, contradicting some aspects and claims that have
been put forward by
some former elaborators on the doctrine/methodology.} Any property of
a physical system
is identified with an {\sl equivalence  class} of {\sl definite experimental
projects}\,\footnote{We deliberately avoid to use the more common
terms like {\it test},
and in particular, {\sl question}, to address definite experimental
projects since
these terms have been the source of confusion and misunderstandings
in some papers (to
which we rather choose not to refer to)\,; a discussion and refutation of these
mistakes can be found in Foulis and Randall (1984).} that can be effectuated on
that system where\,:
\begin{itemize}
\item A definite experimental project is a precisely defined physical procedure
$\alpha$ which includes specification of what should be
conceived as the {\sl positive outcome} when we would effectuate $\alpha$\,;
\item A definite experimental project $\alpha$ is {\sl certain} for a
{\sl particular
realization} of the system, i.e., for the system in a certain {\sl state},
if we obtain the positive outcome {\it with certainty}\, whenever we
{\it would}\, effectuate
$\alpha$ on the system in that particular realization\,;
\item Two definite experimental projects are equivalent whenever
certainty of one is equivalent to
certainty of the other\,; the underlying preorder (or quasi order)
that generates this
equivalence then encodes for two definite experimental projects
$\alpha$ and $\beta$ as
``\,$\alpha\prec\beta$ if and only if certainty of $\alpha$ implies
certainty of $\beta$\,''\,;
the corresponding equivalence class of $\alpha$ will be denoted as
$[\alpha]$ and the physical property
to which it corresponds as $a$\,;
the property $a$ is then called {\sl actual} for a particular
realization of the system, or {\sl true} if one prefers, whenever
$\alpha$ is certain for it.
\end{itemize}
What are the consequences of this operational identification of
properties with definite
experimental projects?  First of all, the partial ordering of the
properties induced by the
preorder on definite experimental projects can now be understood as
an {\it implication
relation} with respect to actuality (or truth).  It also follows that
the properties
constitute a complete lattice $L$\,. Indeed, given $a_i\in L$ with
corresponding
definite experimental projects
$\alpha_i$ we can define the {\sl product} $\prod_i\alpha_i$ as the
definite experimental
project that consists of performing one of the $\alpha_i$, chosen in
{\it any possible way}\,.
The property then defined  by
$[\prod\alpha_i]$ is true if and only if {\it each}\, of the
$\alpha_i$ is true and can
therefore be understood as the {\sl conjunction} of $\{a_i\}_i$\,.
It then also
obviously follows that with respect to the above discussed
preordering of properties,
$[\prod\alpha_i]$ is indeed the {\sl meet} $\bigwedge_i a_i$ of
$\{a_i\}_i$ in the complete lattice $L$\,. By Birkhoff's theorem all
subsets of the lattice then
also have a least upper bound given by
$$
\bigvee_i a_i:=\bigwedge\{b\in L|\forall {i}: a_i\leq b\}\,.
$$
Note here that contra the
usual motivation that conjunctions should be finite, our operational
methodology motivates
arbitrary infinitary ones.  Moreover, as already announced above,
given a meet-semilattice $L$
in which the elements are all the properties of a physical system
ordered by implication with
respect to actuality, and in which the meets encode conjunctions, the products
$\prod\alpha_i$ and corresponding properties
$\bigwedge_i a_i$ then provide an interpretation for the
supplementary elements in the
canonical or MacNeille completion of $L$.\footnote{For an outline of
canonical or
MacNeille completion see Banaschewki and Bruns (1967).}
We will come back to this point, which is slightly more subtle than
it might look at first, in
Section 4.2. Conclusively, although in a considerable number of papers
a property lattice or its
abstract counterpart, somewhat abusively called an {\sl algebraic
quantum logic}, is
conceived as an orthomodular lattice,\footnote{For orthomodular
lattices see Kalmbach (1983) and Bruns and Harding (2000).  An
explicit definition can also be
found in Section 3.2 of this paper.} not necessarily
complete, we will initially consider property
lattices as being general complete lattices. Complete ortholattices
are then a particular
species.

\smallskip
We will now discuss the join in property lattices.
First recall that in the intuitionistic sense, truth of a {\sl disjunction}
coincides with
truth of one of its members, and it is as such that we will conceive
disjunction from now
on. Referring to orthodox Hilbert space quantum mechanics, the
properties of a physical system
are represented by the closed subspaces of a Hilbert space
${\cal H}$ and the join encodes as the closed linear span. So the
join of two atomic
properties $p_1$ and $p_2$ represented by two non-equal
rays\,\footnote{We will somewhat
abusively denote rays in ${\cal H}$ by a representative unit vector.}
$\phi_1$ and
$\phi_2$ in
${\cal H}$ is implied (in the above discussed operational sense) by
any atomic property
$q$ encoded as a ray $\psi$ in the plane spanned by $\phi_1$ and
$\phi_2$\,.  Thus, if
$q\not=p_1,p_2$ then $q\leq p_1\vee p_2$ although actuality of $q$
excludes that of $p_1\!$ {\it and}\, it
excludes that of
$p_2$, i.e., it excludes actuality of
$p_1\!$ {\it or}\,
$p_2$\,.\footnote{Note here that one of the De Morgan laws, namely
$\neg(a\vee b)=\neg
a\wedge\neg b$, indeed still holds in an intuitionistic setting.
However, $\neg(a\wedge
b)=\neg a\vee\neg b$ is not valid anymore.} As such, the join
is in general not a disjunction, and this observation lies at the
base of the construction
made in this paper.  This fact that the join, which (as we saw above)
is from an operational
perspective defined in a secondary way via the meet, is not a
disjunction has been used
by Aerts (1982) to encode the to physicists well-known notion of {\sl
superposition}\,:\,\footnote{In many {\it interpretations}\, of
orthodox quantum theory, a
superposition is understood as a decomposition $\oplus_ic_i\phi_i$ of
a ray $\psi\in{\cal
H}$ that represents the initial state $q$ (so
$\psi=\oplus_ic_i\phi_i$), where $\{\phi_i\}_i$ is an orthonormal
base of ${\cal H}$ that
represents the possible outcome states $\{p_i\}_i$ of a measurement,
envisioning the states
$p_i$ with non-zero $c_i$ as the ``possible truths after the
measurement'' whenever the
measurement will have been effectuated. E.g., Schr\"odinger's cat
ged\"anken experiment\,: The
cat is neither dead nor alive but in a superposition, say {\it
dead}$\oplus${\it alive}, as
long as the measurement is not completely effectuated, where this
effectuation in particular
includes observing whether the cat is dead or alive.  We don't
subscribe to this perspective
but envision a superposition state $q$ of $p_1$ and $p_2$, just as a
different possible
realization of the system,  where  $p_1$, $p_2$ and $q$ are related
by the fact that
$p_1\vee p_2=q\vee p_2=p_1\vee q=p_1\vee p_2\vee q$. In this view the quantum
mechanical measurement described above then induces a change of the
state $q$ to a
state in $\{p_i\}_i$ for which
$c_i\not=0$.} If two states of the system are represented by the atomic
properties
$p_1$ and $p_2$, then all other states $q$ such
that $q\leq p_1\vee p_2$ are called superpositions of $p_1$ and
$p_2$\,. Note here however that
this notion of superposition has been introduced in Aerts (1982)
under the paradigm that states
are indeed in one to one correspondence with atomic properties and
that these atomic properties
are join dense in the property lattice, i.e., the property lattice is
atomistic.  We stress
that we do not  fully subscribe to
this paradigm\,; for a counterexample that employs non-atomistic
property lattices within the
operational methodology we refer to Coecke (2000).
Other examples emerge by restriction of the property lattice led by
certain topological
considerations, for which we refer to Section 4.2 in this paper.
Nevertheless, this notion of
superposition can be extended to non-atomistic property lattices or
more general, any situation
where the states are not encoded as properties, in the following
way\,: Actuality of
$a\vee b$ does not necessarily imply actuality of $a$ or actuality of
$b$\,, i.e., there exists a state for which $a\vee b$ is
actual, but neither $a$ nor
$b$ are actual. We can as such define the following for $A\subseteq L$\,:
\begin{itemize}
\item
{\sl Superposition states} introduced by the join of $A$ are
those states for which $\bigvee A$ is actual while no $a\in A$ is actual.
\end{itemize}
One could as such say that an aspect that characterizes
quantum(-like) property lattices is that {\it the join introduces
superpositions}. This
introduction of superpositions by the join should be conceived as the
strict counterpart of a
join that indeed behaves as a disjunction\,: Actuality of a join does
not necessarily coincide
with actuality of its members since its actuality might be implied by
one of the
superpositions it introduces. Besides superposition states we will
also need to consider the
following for
$A\subseteq L$\,:
\begin{itemize} 
\item
{\sl Superposition properties} introduced by the join of $A$ are
those $c<\bigvee A$ whose
actuality doesn't imply that at least one $a\in A$ is actual, i.e.,
for which there exists a
state that makes $c$ actual while no $a\in A$ is actual.
\end{itemize}
Given the so called {\sl Cartan map} $\mu: L\to{\cal P}(\Sigma)$ with ${\cal P}(\Sigma)$
the powerset of the state set and which assigns to any
property the states in which it is actual, the join of $A$
introduces superposition states if and only if $\mu(\bigvee
A)\not=\bigcup\mu[A]:=\bigcup_{a\in
A}\mu(a)$\,,\footnote{We will use square brackets as a notation for
{\sl pointwise application of
a map} throughout the paper, i.e., $f[X]:=\{f(x)\,|\,x\in X\}$\,.} and
$c<\bigvee A$ is a superposition property if and only if
$\mu(c)\not\subseteq\bigcup\mu[A]$\,.
We will now study some properties of these Cartan maps that will be
of use in this paper.
To distinguish between preservation of finite meets and arbitrary
meets (including the empty
meet $\bigwedge\emptyset=1$\,) we will refer to the latter as an {\it
inf}-morphism.
Since any injective {\it inf}-morphism is also an {\sl order embedding}, i.e.,
$f(a)\leq f(b)\Leftrightarrow a\leq b$\,, we will call it an {\it
inf}-{\sl embedding}.
Call an {\it inf}-morphism $f:L\to L'$ {\sl balanced} if $f(0_L)=0_{L'}$\,.
\bpr
$\mu: L\to{\cal P}(\Sigma)$ is a balanced {\it inf}-embedding.
\epr
\bpf
If $\mu(a)=\mu(b)$, then actuality of $\alpha$ coincides with that of
$\beta$, i.e.,
$\alpha\in[\beta]$ so $a=b$\,. Preservation of infima follows from the
construction of infima
via products, i.e., they stand for conjunction with respect to
actuality. Since the bottom $0$
stands for {\it the absurd}\, it cannot be actual in any state so
$\mu(0)=\emptyset$\,.
Since the top $1$ stands for {\it the trivial}\, it is actual in any state
so $\mu(\bigwedge\emptyset)=\mu(1)=\Sigma=\bigcap\emptyset$\,.
\epf
{\it\underline{Example\,:}}

\smallskip\noindent
Setting $\mu(a):=\{p\in\Sigma|p\leq a\}$ for a complete atomistic
lattice $L$ with atoms
$\Sigma$ we have by atomisticity that $a=\bigvee\mu(a)$\,. Since
$\mu(\bigwedge A)=
{\{p\in\Sigma\,|\,p\leq\bigwedge A\}}=
{\{p\in\Sigma\,|\,\forall a\in A:p\leq a\}}=
\bigcap_{a\in A}{\{p\in\Sigma\,|\,p\leq a\}}=
\bigcap\mu[A]$, since $\mu(0)=\emptyset$ and
since $\mu(a)=\mu(b)\Rightarrow\bigvee\mu(a)=\bigvee\mu(b)\Rightarrow a=b$
it follows that
$\mu:L\to{\cal P}(\Sigma):a\mapsto\mu(a)$ is a balanced {\it inf}-embedding.

\medskip\noindent
Clearly, the Cartan map captures as such the essence of operational
methodology.
In particular are conjunctions now encoded as intersections in the
state space, i.e.,
$\mu(\bigwedge A)=\bigcap\mu[A]$ and disjunctions $\bigvee A$
exactly coincide with
unions, i.e.,
$\mu(\bigvee A)=\bigcup\mu[A]$\,. Note that
$\mu(a)=\emptyset$ implies
$a=0$ by injectivity.
 From injectivity of $\mu$ it also follows that $a<b$ encodes as
$\mu(a)\subset\mu(b)$\,.
Applying all this, superposition properties relate to superposition
states in the following
way\,:
\bpr\label{PropOplus}
If the join of $A\subseteq L$ has superposition properties, then it
also has superposition
states\,; the converse is in general not true.
\epr 
\bpf
 From $c<\bigvee A$ follows $\mu(c)\subset\mu(\bigvee A)$ so if
$\mu(c)\not\subseteq\bigcup\mu[A]$ then $\mu(\bigvee
A)\not\subseteq\bigcup\mu[A]$\, which
proves the first claim. For $L:=\{0,a,a',1\}$ and
$\Sigma:=\{p_1,p_2,q\}$ with $\mu(a)=\{p_1\}$,
$\mu(a')=\{p_2\}$ and
$\mu(1)=\{p_1,p_2,q\}$, the join of $a$ and $a'$ introduces a
superposition state $q$ but no
superposition property.
\epf
{\it\underline{Example\,:}}

\smallskip\noindent
For a complete atomistic lattice $L$ with atoms
$\Sigma$ and $\mu$ as defined in the example above, existence of a
superposition state $p$
does imply existence of a superposition property, namely $p$ itself.

\medskip
Since we want to use the operational methodology as a motivation for
a construction starting
from a complete lattice envisioned as a property lattice, we have to
make an assumption that the
physical essence with respect to superpositions is fully encoded in
the property lattice 
itself and not just in the Cartan map $\mu$\,. This is indeed
necessary, as the counterexample
in the proof of Proposition
\ref{PropOplus} shows\,: Even for the simplest example of a
non-trivial complete Boolean
algebra, namely the square
$\{0,a,a',1\}$, the join $a\vee a'$ is not necessarily interpretable
as a disjunction if
$\mu$  is arbitrarily chosen.
Therefore we will assume at this point that the converse of Proposition
\ref{PropOplus} is also true, i.e., existence of superposition states
implies that of
superposition properties. Denoting the superposition states
introduced by the join of $A$ as
$S_\oplus(A)$ and the superposition properties as $L_\oplus(A)$ this
translates as\,:
\begin{itemize}
\item $p\in S_\oplus(A)\ \Rightarrow\ \exists c_p\in L_\oplus(A):p\in\mu(c_p)$
\end{itemize}
an axiom to which we will refer as {\sl superpositional
faithfullness} of the property lattice (w.r.t.~some Cartan map which is not explicitely
specified ). We will discuss the interpretation and consequences of posing it at the end of
this paper.  We also will
investigate what happens when we drop it, and why the considerations
made in this paper will
then still be usefull.  

\smallskip
Besides this fact that in a property lattice joins do not behave as
disjunctions, the {\it emergence of disjunction in
measurements} is exactly one of the core ingredients of quantum
theory.  Indeed, any
measurement on a system that is not in an eigenstate of that
measurement changes the state of
the system in a non-deterministic manner.  The resulting outcome
state will as such be a
member in a set of possible outcome states. Put in terms of
properties, actuality of a
property $a$ before the measurement guarantees that either
$b_1\!$ {\it or}\, $b_2\!$ {\it or}\, $b_3\!$ {\it or} $\ldots$ will
be actual after the
measurement.  We define an {\sl actuality set} as a set of properties
in which at least one
member is actual. These actuality sets should then be conceived as
the logical propositions
that encode disjunction of actuality of properties, where disjunction
is now indeed to be
understood in the intuitionistic sense, i.e., truth of a disjunction
$A$ ($=$ $A$ is an
actuality set) coincides with actuality ($=$ truth) of one of its
members.  Thus, in other
words, actuality sets recapture the notion ``actuality'' in the
passage from properties to
propositions.  Note that in this setting it is obvious to consider
arbitrary infinitary
disjunctions.  On these actuality sets one can now define an {\sl operational
resolution} (Coecke and Stubbe 1999a,b) as a map that assigns to each
actuality set the
strongest property of which the actuality is implied by the actuality
set.\footnote{In Coecke
and Stubbe (1999a,b) operational resolutions are defined in a
slightly different fashion,
namely as a map ${\cal R}:{\cal P}(\Sigma)\to L$\,, assigning to a set
of states the strongest property that is actual for each of the
realizations in this set. Any
such operational resolution on the states then canonically induces
one on the properties. When substituting
$L$ formally by the isomorphic set of the ${\cal R}$-closed subsets
of $\Sigma$ (every operational
resolution indeed factors in a closure operator ${\cal C}:{\cal
P}(\Sigma)\to{\cal P}(\Sigma)$ and an
isomorphism on its range) then the operational resolution can be seen
as the operation that ``adds all
superposition states''\,.}  Formally, given an actuality set
$A$\,, this property is given by
$\bigvee A$\,. As such, the operational resolution recaptures the
operationally induced logical
structure of the properties on the level of actuality sets.
But what should be considered as the logical
structure of these actuality
sets\,; how do their conjunctions and disjunctions encode? Note for
example that in ${\cal P}(L)$\,, the obvious first candidate to encode actuality sets, 
for properties $a\leq b$ we have that $\{a\}\wedge_{{\cal
P}(L)}\{b\}=\{a\}\cap\{b\}=\emptyset$\,, which
clearly doesn't encode conjunction; thus, ${\cal P}(L)$ is inappropriate as a logic of
actuality sets. We will provide a solution to this in Section 3.1. 

\smallskip 
To conclude this section, if we want to describe a physical system by a
``language'' that is
closed under all
disjunctions of  properties, we formally need to introduce those additional
propositions that express disjunctions of properties and that do not
correspond to a
property in the property lattice.  Thus we want to embed the property
lattice within a larger
propositional lattice whose elements are to be interpreted as
actuality sets whenever the
system is in a state that makes the according proposition true.
The next section will provide the mathematical tool that establishes
this embedding in an
{\it optimal}\, and even {\it universal}\, manner.  

\bigskip\noindent
{\bf 2.2. Bruns-Lakser distributive hulls for complete lattices}

\medskip\noindent
As discussed in the previous section, infima play in the property
lattice a fundamental role
having a direct operational and logical interpretation, respectively
via products of definite
experimental projects and as a conjunction (whereas the join is only
secondary defined) encoded
in terms of the Cartan map as a balanced {\it inf}-embedding.

\smallskip 
Whenever an {\it inf}-embedding is an inclusion, then we call its
domain an {\sl
{\it inf}-subobject} of its codomain.\footnote{We are aware of the potential objections
against this designation, in  particular for category theorists.  We however couldn't think
of anything more suitable.} An {\it inf}-subobject is balanced if the corresponding {\it
inf}-morphism is.  Analogously we define a {\it meet}-{\sl morphism} as a map between
meet-semilattices that preserves finite meets and a {\it sup}-{\sl morphism} as a map
between complete lattices that preserves all suprema, including the empty join
$\bigvee\emptyset=0$\,. Recall that a closure (operator)
${\cal C}:L\to L$ on a complete lattice $L$ is isotone, i.e., $a\leq
b\Rightarrow{\cal
C}(a)\leq{\cal C}(b)$, increasing, i.e., $a\leq{\cal C}(a)$, and
idempotent, i.e., ${\cal C}(a)={\cal
C}^2(a)$. It is {\sl normalized} if moreover
${\cal C}(0)=0$. One then obtains that the range ${\cal C}(L)$ of a closure
operator
${\cal C}$ on a complete lattice $L$ is a {\it inf}-subobject of $L$\,,
which is balanced
whenever the closure is normalized: 
We have for all $a\in A\subseteq{\cal C}(L)$ by isotonicity that ${\cal C}(\bigwedge_L
A)\leq{\cal C}(a)=a$\,; since moreover $\bigwedge_L A\leq{\cal C}(\bigwedge_L A)$, the fact
that
$\bigwedge_L A$ is an infimum forces ${\cal C}(\bigwedge_L A)=\bigwedge_L A$\,, and thus
$\bigwedge_L  A=\bigwedge_{{\cal C}(L)} A$\,.  Conversely, any balanced {\it inf}-subobject
$M$ of
$L$ defines a normalized closure operator
${\cal C}_M:L\to L:a\mapsto\bigwedge\{b\in M|a\leq b\}$\,.  Codomain
restriction of a closure
${\cal C}$ to ${\cal C}(L)$ turns it into a {\it sup}-morphism, so for
$A\subseteq{\cal C}(L)$ we have
$\bigvee_{{\cal C}(L)}A={\cal C}(\bigvee_LA)$.\footnote{See also the appendix at the end of
this paper.}

\smallskip
Given a meet-semilattice $H$, i.e., a poset that admits finite meets,
which is also bounded, i.e., it has $0$ and $1$, then we call it a {\sl Heyting
semialgebra}\,\footnote{Such a Heyting semialgebra exhibits all
structural features of a Heyting
algebra, i.e., a Heyting semialgebra that admits all finite joins.
Since within the context of
this paper (finite) joins defined via a separate operation
(and not canonically
related to meets via Birkhoff's theorem) have no status at all we
choose not to include them in this definition.  For a discussion of
Heyting semialgebras we
refer to Coecke, Moore and Smets (2001c).} if and only if there
exists an additional operation
$(-\Rightarrow -):H\times H\to H$ such that $a\wedge b\leq c$ iff
$a\leq (b\Rightarrow c)$\,.
A Heyting semialgebra which is complete (as a lattice) is called a
{\sl complete Heyting algebra}.
We say that a subset $A$ of a meet-semilattice $L$ has a {\sl
distributive join} if (i) its supremum exists, and (ii) for all $b\in
L$ we have
$b\wedge\bigvee A=\bigvee\{b\wedge a|a\in A\}$. We will abbreviate
this by saying that $\bigvee
A$ is distributive. One then verifies that complete Heyting algebras
are exactly
meet-semilattices in which every subset has a distributive
join.\footnote{For proofs we refer
to  Johnstone (1982) or the appendix at the end of this paper.}
The {\sl Heyting implication} $(-\Rightarrow -):H\times H\to H$ then fixes the
{\sl
Heyting negation} as $\neg(-):=(-\Rightarrow 0)$. Algebraically, this
Heyting negation is a {\sl
pseudo-complementation} since in general it does not satisfy one of
the De Morgan's laws, namely
the {\sl excluded middle law} $\neg a\vee a=1$\,.

\smallskip  
We will now formulate the Bruns-Lakser results.  Recall that given a
category, i.e., a class of objects equipped
with compositionally closed sets of morphisms including identities (e.g., meet-semilattices
with {\it meet}-morphisms or complete lattices with {\it inf}-morphisms), an object $H$ is
{\sl injective} if for every
morphism
$f:L\to H$ with $L$ a subobject of $L'$ (e.g., respectively a {\it
meet}-subobject or {\it
inf}-subobject) there exists a domain extension
$f':L'\to H$.  Given a subobject $L$ of $H$, then we call $H$ an {\sl
essential extension} of
$L$ whenever injectivity of the domain restriction of a morphism
$f:H\to L'$ to $L$ implies
injectivity of $f$ itself. An {\sl injective hull} is then an
essential injective extension.
Note that such an essential injective extension is actually a {\it
minimal} inclusion as a
subobject in an injective object. Indeed, if given an injective hull
$H$ of $L$,
and if $H'$ is another injective object that has $L$ as a subobject,
then, by injectivity of
$H$, the inclusion $L\hookrightarrow H'$ extends to a morphism $H\to
H'$, and since $H'$ is an
essential extension $L$ this map is injective, so $H$ is isomorphic
to a subobject of $H'$.
In particular it also follows that injective hulls are unique up to an
isomorphism. Bruns and
Lakser proved that:
\begin{itemize}
\item Injective meet-semilattices coincide with complete Heyting algebras;
\item Every meet-semilattice has an injective hull.
\end{itemize}
They also provided an implicit and explicit characterization of these
injective hulls of
semilattices\,:
\begin{itemize}
\item A complete Heyting algebra $H$ that has $L$ as a {\it
meet}-subobject is the injective hull of $L$ if and only if\,:
\begin{enumerate}
\item $L$ is {\sl join-dense} in $H$, i.e., for all $a\in H:
a=\bigvee_H\{b\in L|b\leq a\}$\,;
\item If $A\subseteq L$ has a distributive join $\bigvee_LA$ then
$\bigvee_HA=\bigvee_LA$\,.
\end{enumerate}
\item The injective hull of a meet-semilattice $L$ is isomorphic to
its collection of {\sl distributive
ideals} ${\cal DI}(L)$ ordered by inclusion, where a distributive
ideal $A\in{\cal DI}(L)$ is
an {\sl order ideal}, i.e., $a\leq b\in A\Rightarrow a\in A$ and
$A\not=\emptyset$, which is
also closed under existing distributive joins, i.e., if $B\subseteq
A$ has a distributive join
then $\bigvee B\in A$\,; the inclusion of
$L$ in an injective hull $H$ then factors as
$L\cong\{\downarrow\! a\,|\,a\in L\}\hookrightarrow{\cal DI}(L)\cong
H$\,, where
the isomorphic correspondence between the distributive hull $H$ of
$L$ and the distributive
ideals ${\cal DI}(L)$ realizes as\,:
\beqa
& &\ \ \theta\ :H\to{\cal DI}(L):a\mapsto\{b\in L|b\leq a\}\\
& &\theta^{-1}:{\cal DI}(L)\to H:A\mapsto\bigvee_H A
\eeqa
\end{itemize}
Note for the implicit characterization that the second condition forces $L$ to
be a balanced {\it
meet}-subobject.  Indeed, since the join of $\emptyset$ is distributive we have
$0_H=\bigvee_H\emptyset=\bigvee_L\emptyset=0_L$\,. To illustrate the
necessity of this second condition
it suffices to consider $L:=\{0,a,a',1\}$ with $a\vee a'=1$ and
$H:=\{0,a,a',b,1\}$ with
$a,a'< b$\,.  The explicit construction shows us that all the
above can be reformulated for
complete lattices and {\it inf}-morphisms\,:\footnote{It was noted in Stubbe (2000), that
this fact can be seen as a particular incarnation of the Yoneda embedding which always
preserves all limits that happen to exist.} Whenever
$H$ is the injective hull of a complete
lattice $L$ envisioned as a meet-semilattice, then $L$ is an {\it
inf}-subobject of
$H$ since
$\{\downarrow\! a\,|\,a\in L\}\hookrightarrow{\cal DI}(L)$ also
preserves arbitrary infima.
Recalling that the {\sl MacNeille completion} of any poset $L$
consists of closing its
principal ideals
$\{\downarrow\! a\,|\,a\in L\}$ under intersections (Banaschewski and
Bruns 1967), one verifies
that, up to an isomorphism, the inclusion of a meet-semilattice $L$ in its
injective hull factors in (i) the MacNeille completion $\bar{L}$ of
$L$ and (ii) the inclusion
of $\bar{L}$ in its injective hull with respect to complete
lattices and {\it inf}-morphisms. For our purpose this {\it
inf}-restriction suffices.
However, for the sceptici concerning the existence of arbitrary
infima in property lattices we
stress that everything also applies both to meet-semilattices and {\it
meet}-morphisms and to complete
lattices and {\it meet}-morphisms. One then gets for free that the
embedding of a complete
lattice in its injective hull is always an {\it inf}-morphism.
For obvious reasons we will refer to all these equivalent injective
hulls as the {\sl
distributive hull} of a complete lattice. We conclude this section
with an example.
\blm\label{4equiCond}
The following are equivalent for $A\subseteq L$\,:
\par\smallskip\par\noindent
(i) $b\leq\bigvee A\ \Rightarrow\ b=\bigvee_{a\in A}(b\wedge a)$;
\par\smallskip\par\noindent
(ii) $\forall b\in L\,:\,b\wedge\bigvee A=\bigvee_{a\in A}(b\wedge a)$.
\elm
\bpf
(i)$\Rightarrow$(ii): We always have
$b\wedge\bigvee A\geq\bigvee_{a\in A}(b\wedge a)$ and from
$b\wedge\bigvee A\leq\bigvee A$
follows by (i) that $b\wedge\bigvee A=\bigvee_{a\in
A}\bigl((b\wedge\bigvee A)\wedge
a\bigr)\leq\bigvee_{a\in A}(b\wedge a)$\,.
(ii)$\Rightarrow$(i): If $b\leq\bigvee A$, then $b=b\wedge\bigvee A$
so $b=\bigvee_{a\in
A}(b\wedge a)$\,.
\epf
{\it\underline{Example\,:}}

\smallskip\noindent
If $L$ is a complete atomistic lattice with $\Sigma$ as atoms then
$\bigvee\{p\in\Sigma|p\leq a\}$ is distributive.
Indeed, (formally) setting $\mu(a):=\{p\in\Sigma|p\leq a\}$ we have
$a=\bigvee\mu(a)$ and
$a\leq b\Leftrightarrow\mu(a)\subseteq\mu(b)$\,, so
$b\leq\bigvee\mu(a)$ implies
$b=\bigvee\mu(b)=\bigvee_{p\in\mu(b)}(b\wedge
p)=\bigvee_{p\in\mu(a)}(b\wedge p)$\,.
Moreover, for the distributive hull $H$ of $L$ we have $H\cong{\cal
P}(\Sigma)$, i.e.,
the distributive hull of a complete atomistic lattice is a complete
atomistic Boolean
algebra.  Indeed, in terms of
${\cal DI}(L)$, consider
\beqa
& &\ \ \theta\ :{\cal DI}(L)\to {\cal P}(\Sigma):A\mapsto A\cap\Sigma\\
& &\theta^{-1}:{\cal P}(\Sigma)\to{\cal DI}(L):T\mapsto\{a\in
L|\mu(a)\subseteq T\}
\eeqa
Since $\mu(a)\subseteq A$ implies $a=\bigvee\mu(a)\in A$ we have
$a\in A\Leftrightarrow\mu(a)\subseteq
A\Leftrightarrow\mu(a)\subseteq A\cap\Sigma$
and thus
$\theta^{-1}\bigr(\theta(A)\bigl)=\{a\in L|\mu(a)\subseteq
A\cap\Sigma\}=A$\,. From $p\in
T\Leftrightarrow\mu(p)\subseteq T\Leftrightarrow p\in\{a\in
L|\mu(a)\subseteq T\}$ follows
$\theta\bigl(\theta^{-1}(T)\bigr)={\Sigma\cap\{a\in L|\mu(a)\subseteq
T\}}=T$.  Thus
$\theta$ and $\theta^{-1}$ are inverse, and since they are isotone
they define an isomorphism.

\bigskip\noindent
{\bf 3. MAIN RESULTS}

\medskip\noindent
This section constitutes the main argument and constructions.

\bigskip\noindent
{\bf 3.1. Complete lattices and operational resolution}

\medskip\noindent
Let us denote the subsets of $L$ that have a distributive join as
${\cal D}(L)$.
In the next proposition we investigate how the existence of
superpositions relates to
distributivity.
\bpr\label{propDisThenDistr}
For $A\subseteq L$ we have
$$S_\oplus(A)=\emptyset\ \Rightarrow\ A\in{\cal D}(L)\,.$$
Assuming superpositional faithfullness of $L$ we moreover have
$$A\in{\cal D}(L)\ \Rightarrow\ S_\oplus(A)=\emptyset\,.$$
\epr
\bpf
Let $c\leq\bigvee A$. From $c\geq\bigvee_{a\in A}(c\wedge a)$ it follows that
$\mu(c)\supseteq\mu\left(\bigvee_{a\in A}(c\wedge a)\right)$.
Since $p\in\mu(c)$ implies existence of $a\in A$ such that
$p\in\mu(a)$ and thus
$p\in\mu(c\wedge a)$ so $p\in\mu\left(\bigvee_{a\in A}(c\wedge a)\right)$,
it follows that
$\mu(c)\subseteq\mu\left(\bigvee_{a\in A}(c\wedge a)\right)$. By injectivity of
$\mu$ this results in $c=\bigvee_{a\in A}(c\wedge a)$, and by Lemma
\ref{4equiCond}
this completes the proof of the first claim.
For a proof of the second statement see Proposition \ref{6equicond}
in Section 4.1,
$(i)\Rightarrow(iii)$\,.
\epf
Thus, under the assumption of superpositional
faithfulness of $L$ disjunctivity and
distributivity of properties coincides. This justifies the point of view that, using the
Bruns-Lakser results:
\begin{itemize}
\item The inclusion of a property lattice $L$ in its distributive
hull $H$ adds to the property lattice all propositions
that express disjunctions of properties.
Indeed, given $A\subseteq L$,
then $\bigvee_H A$ expresses this disjunction since all suprema in a
complete Heyting algebra
are distributive.
\item It does this in a non-redundant way.  Indeed, by the implicit
characterization it
follows that (i) existing disjunctions are preserved, and (ii) any
other element $a\in H$
indeed expresses a disjunction of properties, namely that of $\{b\in
L\,|\,b\leq a\}$\,.
\item This embedding preserves (i) all infima of
properties, i.e., all conjunctions, (ii) the trivial and (iii) the
absurd, since $L\hookrightarrow H$ is a balanced {\it inf}-embedding.
\end{itemize}
Thus we have embedded $L$ in a logic of propositions that goes equipped with a
pseudo-complementation and internal implication arrow that satisfies
the same rules of
definition and inference of intuitionistic logic.  Moreover, as we
will see below, conjunctivity
and disjunctivity for properties will lift to conjunctivity and
disjunctivity for propositions
such that we indeed have embedded the property lattice in a true
intuitionistic logic which
as such goes equipped with
an intuitionistic negation and implication. However, the inclusion of
$L$ provides $H$ with an
additional operation, namely a normalized closure
\beqa
{\cal R}:H\to H:a\mapsto\bigwedge_H\{b\in L|b\geq a\}\!
&=&\!\bigwedge_L\Bigl\{b\in L\Bigm|b\geq\bigvee_H\{c\in L|c\leq a\}\Bigr\}\\
&=&\!\bigwedge_L\Bigl\{b\in L\Bigm|\forall d\in\{c\in L|c\leq a\}:
b\geq d\Bigr\}\\
&=&\!\bigvee_L\{c\in L|c\leq a\}
\eeqa
referred to as the {\sl operational resolution}\,, an operation
which recuperates the logical structure of properties on the level of
propositions.  

\smallskip 
This operational resolution can indeed be seen as a domain extension and
codomain restriction up to isomorphism of the
operational resolution in Coecke and Stubbe (1999a,b) discussed
above, in the sense that it
assigns the strongest property implied by a proposition, i.e., in
terms of distributive ideals,
implied by an actuality set.
The explicit characterization of distributive ideals will indeed
enable us to envision the above in
terms of actuality sets. First note that $A\subseteq L$\,, as an actuality
set, is equivalent both to
the {\sl implicative closure} $\downarrow\![A]$ of $A$ and {\sl
disjunctive closure}
$\{\bigvee_LB|B\subseteq A\cap{\cal D}(L)\}$ of $A$\,, respectively because of the
implicative significance of the $L$-ordering and disjunctivity of
distributive $L$-suprema. It makes as such sense to consider the distributive ideals ${\cal
DI}(L)$ as the suprema, with respect to inclusion, of equivalence classes of
actuality sets for the following relation\,: Since ${\cal DI}(L)$ is closed under
intersections we can define the closure
$${\cal C}:{\cal P}(L)\to{\cal P}(L):A\mapsto\bigcap\{B\in{\cal
DI}(L)|B\supseteq A\}$$
and an equivalence relation $\sim\,\subseteq\!{\cal P}(L)\times{\cal
P}(L)$ by $A\!\sim\!
B\Leftrightarrow{\cal C}(A)={\cal C}(B)$.
The logical connectives on propositions then translate into a logic of
actuality sets\,:
\beqa
\bigwedge_{{\cal DI}(L)}\!&:&\!{\cal P}({\cal DI}(L))\to {\cal DI}(L):{\cal 
A}\mapsto\bigcap{\cal A}\\
\bigvee_{{\cal DI}(L)}\!&:&\!{\cal P}({\cal DI}(L))\to {\cal
DI}(L):{\cal A}\mapsto
{\cal
C}\!\left(\bigcup{\cal A}\right)\\
\Rightarrow_{{\cal DI}(L)}\!&:&\!{\cal DI}(L)\times{\cal
DI}(L)\to{\cal DI}(L):(B,C)\mapsto
\bigvee_{{\cal DI}(L)}\{A\in{\cal DI}(L)|A\cap B\subseteq C\}\\ 
&&\!\hspace{6.0cm}=\{a\in L|\forall b\in B:a\wedge b\in C\}\\
\neg_{{\cal DI}(L)}\,\!&:&\!{\cal DI}(L)\to{\cal
DI}(L):A\mapsto(A\Rightarrow\downarrow\!0)\\
{\cal R}_{{\cal DI}(L)}\!&:&\!{\cal DI}(L)\to{\cal
DI}(L):A\mapsto\bigm{\downarrow}\!\Bigl(\bigvee_L A\Bigr)
\eeqa 
For ${\cal A}\subseteq{\cal DI}(L)$ we then have that $\bigvee_{{\cal
DI}(L)}{\cal A}={\cal
C}\!\left(\bigcup{\cal A}\right)$ is an
actuality set if and only if at least one $a\in\bigcup{\cal A}$ is
actual, this since all
elements $\bigvee_LB$ in the disjunctive closure are distributive,
and as such if and only if at
least one $A\in{\cal A}$ is an actuality set.  Thus, $\bigvee_{{\cal
DI}(L)}$ is disjunctive, hence $\bigvee_H$ also. If $\bigwedge_{{\cal DI}(L)}{\cal A}$ is an
actuality set then at least
one $a\in\bigcap{\cal A}$ is actual so all $A\in{\cal A}$ are
actuality sets. Conversely, if
all $A\in{\cal A}$ are actuality sets then for all $A\in{\cal A}$ at
least one $a_A\in A$ is
actual such that ${\bigwedge_L\{a_A|A\in{\cal A}\}}$ is actual and
thus ${\bigm\downarrow\!\Bigl(\bigwedge_L\{a_A|A\in{\cal
A}\}\Bigr)}={\bigcap\{\downarrow\!a_A|A\in{\cal A}\}}\subseteq{\bigcap{\cal
A}}={\bigwedge_{{\cal DI}(L)}{\cal A}}$ is an actuality set.  Thus,
$\bigwedge_{{\cal DI}(L)}$
is conjunctive, and as such $\bigwedge_H$ also. This then proves the
claim made above that we
have embedded $L$ in a true intuitionistic logic. 

\smallskip
We are now at the point to understand what an {\sl implication arrow
on properties} should
be.\footnote{We refer to the appendix at the end of this paper for
other attempts which did not
succeed in capturing what we conceive as the true nature of implication.}  It canonically
turns out to be an
external operation
$$
\Rightarrow_L\ :\ L\times L\to{\cal P}(L):(b,c)\mapsto
\{a\in L|a\wedge b\leq c\}
$$
obtained by restricting $\Rightarrow_{{\cal DI}(L)}$\,. If and only
$L$ is itself a complete
Heyting algebra, then we can represent this external operation
faithfully as an internal one by
setting $(b\Rightarrow c):=\bigvee_L(b\Rightarrow_L c)$\,. In
particular our external
implication arrow can be defined by
$$
a\wedge b\leq c\,\Leftrightarrow\, a\in(b\Rightarrow_L c)\,,
$$
as such in a more explicit manner expressing that it generalizes the
implication that lives on
a complete Heyting algebra where $a\in(b\Rightarrow_L c)$ then coincides with
$a\leq(b\Rightarrow c)$\,. The set $(b\Rightarrow_L c)$ is then indeed
the set of properties whose actuality
makes the deduction ``if $b$ is actual then $c$ is actual'' true, and
this is exactly the transcription in
terms of actuality of the minimal requirement of any {\sl functional
formal implication} with respect to
extensional quantification over the state set, i.e., given
$a\in(b\Rightarrow_L c)$\,, then $\forall
p\in\mu(a):p\in\mu(b)$ {\it implies} $p\in\mu(c)$\,.

\smallskip
An {\it inf}-subobject $L$ of $H$ is {\sl distributive join dense} in $H$,
denoted as ${\cal DJD}$\,, if it is join dense in $H$ and if $\forall
a\in H:{\{b\in L|b\leq
a\}}\in{\cal DI}(L)$. A closure ${\cal F}:H\to H$ is ${\cal DJD}$ if
${\cal F}(H)$ is ${\cal
DJD}$ as an {\it inf}-subobject of $H$.
Note that since $\emptyset\not\in{\cal DI}(L)$ the {\it
inf}-subobject inclusion is balanced so
a ${\cal DJD}$-closure is always normalized. Referring back to the
implicit Bruns-Lakser
characterization of distributive hulls  of the previous section, the
requirement $\forall a\in
H:{\{b\in L|b\leq a\}}\in{\cal DI}(L)$ is equivalent to the inclusion ${\cal
F}(H)\hookrightarrow H$ preserving existing distributive joins.
Indeed, given that for
$A\in{\cal D}\left({\cal F}(H)\right)$ we have $\bigvee_{{\cal
F}(H)}A=\bigvee_HA$, or equivalently
$\bigvee_{{\cal F}(H)}A\leq\bigvee_HA$, then $A\subseteq\{b\in {\cal
F}(H)|b\leq
a\}$ implies $\bigvee_{{\cal F}(H)}A\leq\bigvee_HA\leq a$\,. Conversely, since
$\{b\in {\cal F}(H)|b\leq\bigvee_HA\}\in{\cal DI}\bigl({\cal
F}(H)\bigr)$ implies
$\bigvee_{{\cal F}(H)}A\in{\{b\in {\cal F}(H)|b\leq\bigvee_HA\}}$ it
follows that $\bigvee_{{\cal
F}(H)}A\leq\bigvee_HA$\,.
We are now in a position to summarize the above within the following
definition\,:
\bdf\label{TH:1}
By the ``intuitionistic or disjunctive representation of quantum
logic'' we refer to
bijective correspondence between isomorphism classes of

\smallskip\noindent
(i) complete lattices, denoted {\sf CLat}\,, and,

\smallskip\noindent
(ii) complete Heyting algebras equipped with a ${\cal DJD}$-closure,
denoted {\sf DJDHeyt}\,,

\smallskip\noindent
which is realized by the following equivalence\,:\,\footnote{An isomorphism
between complete Heyting algebras $H_1$ and $H_2$ equipped with
respective ${\cal
DJD}$-closures ${\cal F}_1$ and
${\cal F}_2$ is an order-isomorphism $h:H_1\to H_2$ such that ${\cal F}_2\circ
h=h\circ{\cal F}_1$.}
\beqa
\theta\ &:&{\sf CLat}\to{\sf DJDHeyt}:L\mapsto\left({\cal
DI}(L),{\cal R}_{{\cal DI}(L)}\right)\\
\theta^*&:&{\sf DJDHeyt}\to{\sf CLat}:(H,{\cal F})\mapsto{\cal F}(H)
\eeqa
Given a complete lattice of ``properties'', the ${\cal DJD}$-closure
operator that arises on
the complete Heyting algebra of ``propositions'', is called
the ``operational resolution''.  It assigns to a proposition the
strongest property
implied by it.
\edf
{\it\underline{Example\,:}}

\smallskip\noindent
Isomorphism classes of complete atomistic lattices and
complete atomic Boolean algebras equipped with a ${\cal DJD}$ closure
are in bijective
correspondence.
Indeed, if $L$ is atomistic then
${\cal DI}(L)\cong{\cal P}(\Sigma)$ and by the Lindenbaum-Tarski theorem these
are exactly complete atomistic Boolean algebras.  Conversely, if
$H\cong {\cal P}(\Sigma)$
then
${\cal DJD}$ requires $p\in{\cal F}(H)$, so ${\cal F}$ is a $T_1$
closure, i.e., all points are
closed, and thus
${\cal F}(H)$ is atomistic via the bijective correspondence of isomorphism
classes of $T_1$-closure spaces and complete atomistic lattices.
Note here that both for orthodox Hilbert space quantum mechanics and
phase space classical
mechanics atomisticity is an axiom so this example covers essentially
the primitive\,\footnote{By primitive
we refer to the fact that for topological, probabilistic or other
reasons one might consider restrictions
of this primitive complete atomistic setting that are not complete or
not atomistic anymore, as we will
discuss in Section 4.2 of this paper.} situations presently encountered
in {\it orthodox physical
theories}\,.  However, since atomisticity cannot be motivated within
the operational methodology as it is
applied in this paper, it shouldn't play a role in any derivation or
construction, and may only be
injected at the end as a particular feature of these paradigm
examples whenever one wants to consider them
explicitly.  Moreover, as we mentioned above and will discuss below, there are indeed
situations where atomisticity is not the case.

\bigskip\noindent
{\bf 3.2 Complete ortholattices and operational complementation.}

\medskip\noindent
We will now go back to our initial goal of merging disjunctive suprema and
orthocomplementation within one structure, as such eliminating the
Birkhoff-von Neumann dilemma.
Note that at this point we do not attribute a particular operational
or physical significance
to this orthocomplementation.  We will just {\it refine}\, the
results above in case that
there is a given one, as it is the case for orthodox quantum theory\,.

\smallskip
A complete lattice $L$ is a {\sl complete ortholattice} if it goes
equipped with an
{\sl orthocomplementation}\ \ $':L\to L$\,, i.e., an operation that
satisfies\,:
\beqa
\mbox{OC1:}\ \ a\wedge a'=0,
\quad\quad
\mbox{OC2:}\ \ a\leq a'',
\quad\quad
\mbox{OC3$_l$:}\ \ a\leq b\Leftarrow b'\leq a'.
\eeqa
It is a {\sl complete pseudo-ortholattice} if $\ ':L\to L$
satisfies OC1, OC2 and
\beqa
\mbox{OC3$_r$:}\ \ a\leq b\Rightarrow b'\leq a'.
\eeqa
It is a {\sl complete ${\cal
DJD}$-pseudo-ortholattice} if in addition the range of the
orthocomplementation is ${\cal DJD}$.
Isomorphisms of complete ortholattices and complete
pseudo-ortholattices are then obviously
those order-isomorphisms that preserve the orthocomplementation.
\bpr\label{axiomequi}
We have the following for the above axioms:

\smallskip\noindent
(i) $[\, $OC2\, , OC3$_l\, ]$ is equivalent to $[\, a=a''$ ,\, $a\leq
b\Leftrightarrow b'\leq a'\, ]$,

\smallskip\noindent
(ii) $[\, $OC2\, , OC3$_r\, ]$ implies $[\, a'=a'''$ ,\, $a'\leq
b'\Leftrightarrow b''\leq a''\, ]$,

\smallskip\noindent
(iii) $[\, a'\wedge a''=0$\, , OC2$\, ]$ implies OC1.
\epr
\bpf
(i): From $a'\leq a'''$ follows $a\geq a''$ by OC3$_l$; $a\leq
b\Rightarrow a''\leq
b''\Rightarrow b'\leq a'$. (ii): $a'\leq b'\Leftarrow b''\leq a''$ by
OC2 and OC3$_r$, the
rest by applying (i) to
$a'$ and $b'$. (iii): $a\wedge a'\leq a''\wedge a'$ by OC2, so
$a\wedge a'\leq 0$\,.
\epf
\bth\label{orthotheorem}
Isomorphism classes of

\smallskip\noindent
(i) complete ortholattices, denoted ${\sf COLat}$\,, and,

\smallskip\noindent
(ii) complete ${\cal DJD}$-pseudo-ortho Heyting algebras, denoted
{\sf DJDOHeyt}\,,

\smallskip\noindent
are in bijective correspondence via the equivalence
\beqa
\theta\ &:&{\sf COLat}\to{\sf DJDOHeyt}:(L,\, ')\mapsto({\cal DI}(L),^\perp)\\
\theta^*&:&{\sf DJDOHeyt}\to{\sf COLat}:(H,\, ')\mapsto(H',\, '\lfloor_{H'})
\eeqa
where the ``operational complementation'' is defined as
$$
^\perp:{\cal DI}(L)\to{\cal DI}(L):A\mapsto\bigm\downarrow\!(\bigvee_L A)'\,.
$$
In terms of inclusion of $L$ into its distributive hull $H$ this translates as
$$
^\perp:H\to H:a\mapsto {\cal R}(a)'\,,
$$
so the operational complementation has the operational resolution as
its square, establishing
the operational complementation as a refinement of the latter.
\eth
\bpf
Consider $(H,\, '\,)\in {\sf DJDOHeyt}$\,. Then $a''=(a'')''$ via
Proposition \ref{axiomequi} (ii) and $a\leq b\Rightarrow a''\leq b''$
via twice OC3$_r$
assure\ \
$'':H\to H$ to be a closure. Moreover, since $a'=(a')''\in H''$ and
$a''=(a')'$ it
follows that $H'=H''$ so $H'$ is a complete lattice. By Proposition
\ref{axiomequi} (ii) we also have
$(a')\leq (b')\Leftrightarrow (b')'\leq (a')'$, so the domain restriction
$'\lfloor_{H'}$ of $'$ to $H'$ defines an orthocomplementation on
$H'$, and thus $\theta^{-1}$
is well defined.  By $H'=H''$ it also follows that $''$ is a ${\cal
DJD}$-closure on $H$ which
we will denote by ${\cal F}$ --- ref$.$ Definition \ref{TH:1}.
For $(L,\, '\,)\in{\sf COLat}$ we have that
$A^{\perp\perp}=\bigm\downarrow\!\Bigl(\bigvee_L\bigl(\downarrow\!(\bigvee_L
A)'\bigr)\Bigr)'=
{\bigm\downarrow\!(\bigvee_L A)''}={\cal R}_{{\cal DI}(L)}(A)$ so
$^{\perp\perp}$
is a closure and thus $A\subseteq A^{\perp\perp}$\,. Since moreover
$\bigvee_L(-)$ and
$\downarrow\!(-)$ are isotone and $'$ is antitone
$A\subseteq B\Rightarrow
B^\perp\subseteq A^\perp$.  Thirdly,
$A^\perp\wedge B^{\perp\perp}=
{\bigm\downarrow\!(\bigvee_L A)'\wedge\bigm\downarrow\!(\bigvee_L A)}=
{\bigm\downarrow\!\Bigl((\bigvee_L A)'\wedge(\bigvee_L A)}\Bigr)=
\downarrow\!0\,$ by OC1 for $L$, so
$^\perp$ is OC1 on ${\cal DI}(L)^\perp$ and thus by Proposition
\ref{axiomequi} (iii) on ${\cal DI}(L)$, assuring
$^\perp$ to be a pseudo-orthocomplementation. Since
$^{\perp\perp}={\cal R}_{{\cal DI}(L)}$ we have ${\cal
DI}(L)^\perp\supseteq\downarrow\![L]$
so the inclusion of the range of $^\perp$ in ${\downarrow\!(-)}$ forces ${\cal
DI}(L)^\perp=\downarrow\![L]$. Thus $^\perp$ is ${\cal DJD}$ so
$\theta$ is well defined.  By
the above it then follows that
$\theta^{-1}\bigl(\theta(L,\, '\,)\bigr)=\bigl(\downarrow\![L]\,,\,
^\perp\lfloor_{\downarrow\![L]}\bigr)\cong(L,\,'\,)$\,.   One also
straightforwardly
verifies that
$\theta\bigl(\theta^{-1}(H,\, '\,)\bigr)={\Bigl({\cal DI}\bigl({\cal
F}(L)\bigr)\,,\bigm\downarrow_{{\cal F}(H)}\!\bigl(\bigvee_{{\cal F}(H)}
A\bigr)'\Bigr)}\cong(H,\, '\,)$ what completes the proof.
\epf
Note that in terms of the {\sl orthogonality relation} $a\perp
b\Leftrightarrow a\leq b'$
induced by the orthocomplementation we have
$A^\perp=\{b\in L\,|\,\forall a\in A:a\perp b\}$\,. Although in the
strict operational
methodology it is possible to motivate the existence of an
orthogonality relation, the fact
that every property in the property lattice can be written as the
supremum of a biorthogonally
closed subset of the lattice, i.e., the orthogonality relation realizes an
orthocomplementation, is taken as an axiom (Jauch and Piron 1969,
Aerts 1982, Moore 1999).
One could wonder whether orthocomplementation, and consequently, ${\cal
DJD}$-pseudo-orthocomplementation of the distributive hull, can be
obtained in a canonical
manner, without having to assume it.
\bqu\label{quortho}
Given a complete lattice $L$ equipped with an orthogonality relation
$\perp\,\subseteq
L\times L$\,, does there exist an elegant characterization of a
minimal extension of
${(L,\perp)}$ as a
${\cal DJD}$-pseudo-ortho Heyting algebras sensu the role in this
paper of distributive hulls
in the category of complete lattices, and can this be translated in terms of a
faithful representation sensu Definition
\ref{TH:1} and Theorem \ref{orthotheorem}?
\equ

\bigskip\noindent
{\bf 4. FURTHER ANALYSIS AND OPEN PROBLEMS}

\medskip\noindent
We discuss remaining loose ends and further research.

\bigskip\noindent 
{\bf 4.1. Characterization of disjunctivity}

\medskip\noindent
We will now characterize the correspondence between disjunctivity and
distributivity in the
strict operational methodology, i.e., when a Cartan map
$\mu:L\to{\cal P}(\Sigma)$ is
explicitly given.  Set
$S(p):=\bigwedge\{a\in L|p\in\mu(a)\}$ and write
$a<\hspace{-2.2mm}\cdot\ b$ if $a<b$ and $c<b\Rightarrow c\leq a$ for all
$c$\,.
\bpr\label{6equicond}
The following are equivalent:
\par\smallskip\par\noindent
(i) $p\in S_\oplus(A)\ \Rightarrow\ \exists c_p\in L_\oplus(A):p\in\mu(c_p)$\,;
\par\smallskip\par\noindent
(ii) $L_\oplus(A)=\emptyset\ \Rightarrow\ S_\oplus(A)=\emptyset$\,;
\par\smallskip\par\noindent
(iii) $A\in{\cal D}(L)\ \Rightarrow\ S_\oplus(A)=\emptyset$\,;
\par\smallskip\par\noindent
(iv) $A\in{\cal D}(L)\ \Rightarrow\ \left[S(p)=\bigvee A\,\Rightarrow\,p\not\in
S_\oplus(A)\right]$\,;
\par\smallskip\par\noindent
(v) $S(p)=\bigvee A\,\Rightarrow\,S(p)\in A$\,;
\par\smallskip\par\noindent
(vi) For $p\in\Sigma$ we either have\,:
$\left\{\begin{array}{l}
1. \downarrow\!S(p)\setminus\{S(p)\}=\downarrow\!a\ {\rm for}\
0<a<\hspace{-2.2mm}\cdot\ S(p)\,;\vspace{1mm}\\
2.\ S(p)\ {\rm is\ an\ atom\ of}\ L\,.
\end{array}
\right.$
\epr
\bpf
We proceed by proving
(v)$\Rightarrow$(iv)$\Rightarrow$(iii)$\Rightarrow$(i)$\Rightarrow$(ii)$\Rightarrow$(v)
and
(v)$\Leftrightarrow$(vi), where (v)$\Rightarrow$(iv) and
(i)$\Rightarrow$(ii) are both trivial.

(iv)$\Rightarrow$(iii): Let $A\in{\cal D}(L)$ and $p\in S_\oplus(A)$\,. Since
$S(p)\leq\bigvee A$ it follows that
$S(p)=S(p)\wedge\bigvee A=\bigvee_{a\in A}\bigl(S(p)\wedge a\bigr)$\,.
Next, since
$b\wedge\bigvee_{a\in A}\bigl(S(p)\wedge a\bigr)=b\wedge S(p)\wedge\bigvee A=
\bigvee_{a\in A}\bigl(b\wedge S(p)\wedge a\bigr)$ we have
$\{S(p)\wedge a|a\in A\}\in{\cal
D}(L)$\,. However, $p\not\in\mu(a)$ for $a\in A$, so
$p\not\in\mu\bigl(S(p)\wedge a\bigr)$,
what results in $p\in
S_\oplus\bigl(\{S(p)\wedge a|a\in A\}\bigr)$ and this conflicts with (iv)\,.

(iii)$\Rightarrow$(i): There are two possibilities for $p\in S_\oplus(A)$\,:
1. $S(p)=\bigvee A$\,: Since $p\in S_\oplus(A)$ we have $S(p)\not\in A$ so
$A\subseteq{\downarrow\!S(p)}\setminus\{S(p)\}$ and thus
$S(p)={\bigvee\bigl({\downarrow\!S(p)}\setminus\{S(p)\}\bigr)}$\,.
We claim that ${\downarrow\!S(p)}\setminus\{S(p)\}\in{\cal D}(L)$.
Indeed, following Lemma
\ref{4equiCond}, it suffices that $c\leq S(p)$ implies
$c=c\wedge\bigvee\bigl({\downarrow\!S(p)}\setminus\{S(p)\}\bigr)=\bigvee\{{c\wedge
a}|a<S(p)\}$ what is the case. By (iii) we obtain
$S_\oplus\bigl({\downarrow\!S(p)}\setminus\{S(p)\}\bigr)=\emptyset$
what contradicts with
$p\in\mu\bigl(S(p)\bigr)$ and $p\not\in\mu(a)$ for $a<S(p)$\,.
2. $0<S(p)<\bigvee A$\,:
It suffices to set $c_p:=S(p)\in L_\oplus(A)$ since
$p\in\mu\bigl(S(p)\bigr)$\,.

(ii)$\Rightarrow$(v): If for some $p\in\Sigma$ we have
$\bigvee\bigl(\downarrow\!S(p)\setminus\{S(p)\}\bigr)=S(p)$ then
$p\in S_\oplus\bigl(\downarrow\!S(p)\setminus\{S(p)\}\bigr)\not=\emptyset$\,.
However,
$L_\oplus\bigl(\downarrow\!S(p)\setminus\{S(p)\}\bigr)=\emptyset$
since
$a<{\bigvee\bigl(\downarrow\!S(p)\setminus\{S(p)\}\bigr)}$ implies
$a\in{\downarrow\!S(p)}\setminus\{S(p)\}$\,, so
$\bigvee\bigl({\downarrow\!S(p)}\setminus\{S(p)\}\bigr)<S(p)$ by (ii)\,.  Thus
$A\subseteq{\downarrow\!S(p)}\setminus\{S(p)\}$ implies $\bigvee
A<S(p)$\,, so $\bigvee
A=S(p)$ implies $A\not\subseteq{\downarrow\!S(p)}\setminus\{S(p)\}$\,. Since
$A\subseteq\downarrow\!S(p)$ we obtain $S(p)\in A$\,.

(v)$\Leftrightarrow$(vi) Above we proved
$\bigvee\bigl({\downarrow\!S(p)}\setminus\{S(p)\}\bigr)<S(p)$.
It thus follows that
${\downarrow\!S(p)}\setminus\{S(p)\}=
\downarrow\!\bigvee\bigl({\downarrow\!S(p)}\setminus\{S(p)\}\bigr)$ where
$\bigvee\bigl({\downarrow\!S(p)}\setminus\{S(p)\}\bigr)<\hspace{-2.2mm}\cdot\
S(p)$\,.
The converse is easily verified.
\epf
The first of these conditions is what we defined as superpositional
faithfullness, in the
sense of\,: ``the property
lattice fully reflects the systems behavior in terms of superpositions''.
Actually, this condition is implicit in everything that has been done up to
date in quantum logic since to the current authors' knowledge no
construction that explicitly
uses distinct state sets and property lattices as primitive objects have
been considered
(except then for states
being measures on the property lattice, but then the concept of state
is not primitive).  In
the operational methodology one initially takes this into account,
but then ``kills'' the
distinction with the axiom that states encode as a join dense set of
atoms of the property
lattice.  The first explicit constructions probably are those  that
can be found in Coecke
and Stubbe (1999a,b) and Coecke (2000).
The second
condition shows that this superpositional faithfullness can be
formulated in a slightly weaker
fashion. The third condition identifies distributivity and
disjunctivity, and condition four
to six constitute a stepwise characterization of the above in terms
of the properties
$\{S(p)|p\in\Sigma\}$\,. We already mentioned that when for an
atomistic lattice the atoms
are envisioned as states we do have superpositional faithfullness.
Since then $S(p):=p$\,,
this corresponds in the above proposition with the case where (vi).1
is excluded in (vi). We
indeed have the following.\footnote{See Moore (1999) for more details on this and related
matters.} 
\bpr
Given a Cartan map $\mu:L\to{\cal P}(\Sigma)$, the set
$\{S(p)|p\in\Sigma\}$ is join dense in
$L$\,, i.e., $a=\bigvee\{S(p)\,|\,p\in\mu(a)\}$ for all $a\in L$\,.
\epr
\bpf
We have $S(p)\leq
a\Leftrightarrow\bigcap\{\mu(c)|c\in L,
p\in\mu(c)\}={\mu\Bigl(\bigwedge\{c\in
L|p\in\mu(c)\}\Bigr)}=\mu\bigl(S(p)\bigr)\subseteq\mu(a)\Leftrightarrow
p\in\mu(a)$ since
$\mu$ is an injective {\it inf}-morphism. Thus we have
$a\geq{\bigvee\{S(p)\,|\,p\in\mu(a)\}}$\,, and less or equal
saturates into an equality
since
$a>b={\bigvee\{S(p)\,|\,p\in\mu(a)\}}$ both implies $\mu(a)\supset\mu(b)$ and
$p\in\mu(a)\Rightarrow b\geq S(p)\Rightarrow p\in\mu(b)$\,, i.e.,
$\mu(a)\subseteq\mu(b)$\,.
\epf
An example radically different from the atomistic one is a completely
ordered set $L$\,, where we
set $\Sigma:=L\setminus\{0\}$ and $\mu:L\to{\cal
P}(L\setminus\{0\}):0\mapsto\emptyset\,;a(\not=0)\mapsto]0,a]$\,. 

\smallskip
It is clear
that not all complete lattices admit a realization as a property
lattice equipped with Cartan
map that is superpositionally faithful\,.  However, below we will
motivate that in view of
certain topological considerations a much larger class of complete
lattices than the one that
one might expect from the results above admits a meaningful
distributive hull in the sense of
disjunctive completion. Still, even within the setting of this
section the distributive
hull of any complete lattice provides a ``lower bound'' for the {\sl
disjunctive hull}, there
where an obvious ``upper bound'' is {\sl downset completion}\,, i.e., any
extension $\bar{H}$ of $L$ isomorphic to
$$
{\cal I}(L):={\left\{\downarrow\![A]\bigm|A\subseteq L\right\}}
$$
that makes the following diagram commute
\[
\begin{array}{ccc}
L&\hookrightarrow&\bar{H}\\
\hspace{-2mm}{\scriptstyle\downarrow(-)}\!\searrow\hspace{-5mm}&&\hspace{-5mm}\nearrow{\scriptstyle\cong}\\
&{\cal I}(L)
\end{array}
\]
Indeed, within the
context of the strict operational methodology it makes sense to
investigate what
characterizes the disjunctive hull given an arbitrary $\mu:L\to{\cal
P}(\Sigma)$
which not necessarily satisfies superpositional faithfullness.   As it is the
case for ${\cal DI}(L)$\,, ${\cal I}(L)$ is a complete Heyting algebra being
closed under unions and intersections and as such inheriting
distributivity from
${\cal P}(L)$\,. Since disjunctions are in bijective correspondence
with unions of
$\mu[L]$ elements it is clear that the {\sl disjunctive hull} is in general any
extension $H_\mu$ of $L$ isomorphic to
$$
{\cal D}_\mu(L):=\left\{\,\bigcup\mu[A]\Bigm|A\subseteq L\,\right\}\,.
$$
that makes the following diagram commute
\[
\begin{array}{ccc}
L\ \!&\hookrightarrow\ &H_\mu\\
\hspace{-2mm}{\scriptstyle\mu}\!\!\searrow\!\hspace{-5mm}&&\hspace{-5mm}\nearrow
{\scriptstyle\cong}\\
&{\cal D}_\mu(L)
\end{array}
\]
\bpr\label{DILas BOUND}
For every Cartan map $\mu:L\to{\cal P}(\Sigma)$ there exist two balanced {\it
inf}-embeddings $\varphi_\mu:{\cal DI}(L)\to{\cal D}_\mu(L)$ and
$\varepsilon_\mu:{\cal
D}_\mu(L)\to{\cal I}(L)$\,. Moreover,
given a complete lattice
$L$ there exists a Cartan map
$\mu:L\to{\cal P}(\Sigma)$ that realizes ${\cal I}(L)$ as disjunctive
hull, i.e., such that
${\cal D}_\mu(L)\cong{\cal I}(L)$\,. However, this is in general not
the case for
${\cal DI}(L)$\,.
\epr
\bpf
First, note that ${\cal D}_\mu(L)$ is a complete Heyting algebra. Indeed,
we have $\bigcup_{i\in I}(\bigcup\mu[A_i])=\bigcup\mu[\bigcup_{i\in I}A_i]$ and
by complete distributivity\footnote{See for example Johnstone (1982) \S VII
p.278--279.} of
${\cal P}(\Sigma)$ we moreover have $\bigcap_{i\in I}(\bigcup\mu[A_i])=
\bigcup_{(x_i)_i\in X}\left(\bigcap_{i\in I}\mu(x_i)\right)=
\bigcup_{(x_i)_i\in X}\mu(\bigwedge_{i\in I}x_i)$ where
$X=\{(a_i)_i\,|\, a_i\in
A_i\}$\,. Thus, it follows that ${\cal D}_\mu(L)$ is closed under all
unions and
intersections and as such inherits distributivity from ${\cal P}(\Sigma)$\,.
As such, existence of $\varphi_\mu$ is guaranteed since ${\cal
DI}(L)$ is a distributive hull of
$L$ and the codomain restriction of the Cartan map $\mu:L\to{\cal
P}(\Sigma)$ to
${\cal D}_\mu(L)$ defines an {\it inf}-inclusion of $L$ in an injective object.
Next, set $\varepsilon_\mu:T\mapsto\{a\in L\,|\,\mu(a)\subseteq T\}$\,.
Since $a\in\varepsilon_\mu(T)\Leftrightarrow\mu(a)\subseteq T$ it follows that
$a\in\bigcap_{T\in{\cal T}}\varepsilon_\mu(T)\Leftrightarrow\forall T\in{\cal
T}:a\in\varepsilon_\mu(T)\Leftrightarrow\forall T\in{\cal T}:\mu(a)\subseteq
T\Leftrightarrow\mu(a)\subseteq\bigcap{\cal T}\Leftrightarrow
a\in\varepsilon_\mu(\bigcap{\cal T})\,$,
so this inclusion $\varepsilon_\mu$ preserves intersections, i.e., infima.
Given an arbitrary complete lattice $L$, setting $\mu:L\to{\cal
P}(L\setminus\{0\}):0\mapsto\emptyset\,;a(\not=0)\mapsto]0,a]$ we clearly
realize ${\cal I}(L)\cong{\cal D}_\mu(L)$ via $A\leftrightarrow
A\setminus\{0\}\,;\{0\}\leftrightarrow\emptyset$\,. However, the
lattice of open
sets (with respect to the standard topology) of the unit interval
cannot realize
${\cal DI}(L)$\,. Indeed, there are no candidates in this lattice to play the
role of $S(p)$ in view of condition (vi) of Proposition \ref{6equicond}\,.
\epf
\bqu\label{questionM(L)}
Is there some categorical property that elegantly characterizes
${\cal D}_\mu(L)$ in some
category with as objects Cartan maps, sensu the role in this paper of
distributive hulls in
the category of complete lattices? How do the different ${\cal
D}_\mu(L)$ relate for fixed
$L$ and what is the status of ${\cal DI}(L)$ and ${\cal I}(L)$ for
this collection/category?
How do the results of Paseka (1994) on covers in generalized frames
fit in this picture?
\equ
Note that in respect of the second question one can verify that two
canonical choices for
morphisms between ${\cal D}_\mu(L)$ and ${\cal D}_{\mu'}(L)$ present
themself, namely a {\it
inf}-morphism $f:{\cal D}_\mu(L)\to{\cal
D}_{\mu'}(L):T\mapsto\bigcup\{\mu'(a)\,|\,\mu(a)\subseteq T\}$\,, and
a {\it sup}-morphism
$g':{\cal D}_{\mu'}(L)\to{\cal
D}_\mu(L):T\mapsto\bigcup\bigl\{\mu\bigl(S'(p)\bigr)\,|\,p\in
T\bigr\}$\,, which prove to
be adjointly related\,.
\bqu
In the above, and in particular in the proof of Proposition
\ref{DILas BOUND}, it
seems that for ${\cal DI}(L)$ there is a strong connection between complete
distributivity and superpositional faithfullness with respect to some Cartan
map.  Can this be put in a simple picture and could this provide a
simplification of the presentation compared to the one in this paper?
\equ
In particular complete distributivity seems to arise when at the starting
point of the construction we restrict to injective hulls of {\it inf}-lattices.
This setting however requires that all infima are conjunctive,
an assumption that we
will drop in the next section.

\bigskip\noindent 
{\bf 4.2. Finitely conjunctive infima}

\medskip\noindent
Although in the strict sense of the operational methodology outlined earlier
$S(p)=0$ is excluded, this since $\bigwedge\{a\in L|p\in\mu(a)\}$ is a property
that is actual in state $p$ since all
$a$ with $p\in\mu(a)$ are actual in $p$, it does seem to make sense to
consider property lattices where only finite meets are conjunctive in view of
certain topological motivations, even within an operational setting, as such
allowing $S(p)=0$ whenever $S(p)$ is the infimum of an infintary set. This
finitely conjunctive property lattice should then be envisioned as a
restriction
of the {\it true} property lattice. There are for example arguments in terms of
{\sl affirmation} vs.~{\sl refutability} motivating that so called
{\sl finitely
observational properties} are restricted to open sets of states, thus proposing
frames as the corresponding property lattices (Vickers 1989). Consider for
example
$$
\Sigma:=[0,1]\,,\,\bigl(L:=\{T\subseteq[0,1]|T\ {\rm is\ open}\}
\,,\subseteq\bigr)\,,\,
\mu(T):=T\,.
$$
Since all suprema are unions they are all disjunctive.  One should
then envision this lattice as a restriction of ${\cal P}([0,1])$ where now only
finite infima are to be seen as conjunctions contrary to ${\cal
P}([0,1])$ itself
where all infima are conjunctions.  Consequently, the map $\mu:L\to{\cal
P}(\Sigma)$ that assigns to properties the states in which they are
actual is now
a balanced {\it meet}-embedding that also preserves the top, but
which is not necessarily
an {\it inf}-embedding anymore. We will refer to these maps as {\sl weak Cartan
maps}\,.  In principle, the domain of such a weak Cartan map should not even be
a complete lattice, but only a bounded meet-semilattice.
\bpr
Every weak Cartan map $\mu:L\to{\cal P}(\Sigma)$ admits a conjunction
preserving
extension as a Cartan map, namely
$$
\bar{\mu}:\bar{L}_\mu\to{\cal P}(\Sigma)
$$
where
$$
\bar{L}_\mu\cong{\cal C}_\mu(L):=\left\{\,\bigcap\mu[A]\Bigm|A\subseteq
L\,\right\}\,.
$$
is restricted by commutation of
\[
\begin{array}{ccc}
L\!\ &\hookrightarrow\ &\bar{L}_\mu\\
\hspace{-2mm}{\scriptstyle\mu}\!\!\searrow\!\hspace{-5mm}&&\hspace{-5mm}\nearrow
{\scriptstyle\cong}\\
&{\cal C}_\mu(L)
\end{array}
\]
\epr
\bpf
Straightforward verification.
\epf
Since the inclusion $L\hookrightarrow\bar{L}_\mu$ is a completion it
factors over MacNeille completion $L\hookrightarrow\bar{L}$
(Banaschewski and Bruns 1967)\,,
where
$$
\bar{L}\cong\left\{\,\bigcap\downarrow\![A]\Bigm|A\subseteq L\,\right\}\,,
$$
again with the obvious commutation property.
In general however, $\bar{L}$ and
$\bar{L}_\mu$ do not coincide\,: take as a counterexample the
standard topology on an interval
with as states the points of the interval. Thus, completeness does
not imply conjunctivity,
although the converse is true.\footnote{Note that via this
conjunctive completion we obtain for
any weak Cartan map as such
$
{\cal R}:{\cal
P}(\Sigma)\to\bar{L}_\mu:T\mapsto\bigcap\{A\in\mu[L]\,|\,T\subseteq
A\}
$
as the operational resolution sensu Coecke and Stubbe (1999a,b)\,.
Deriving an operational
resolution from a Cartan map indeed requires conjunctivity of all
infima in the property lattice.}
One could say that the distributive hull plays the same role for
disjunctive completion as
MacNeille completion plays for conjunctive completion.
By the above it also follows that it makes no essential difference to
work either with
finitary conjunctive meets or infinitary conjunctive infima.\footnote{Note
that the finitary
representation expresses the  non-primitive nature of suprema in a
much stronger sense since in general they even don't exist.}
When evaluating superpositional faithfullness for property lattices
with finitely conjunctive
infima the conditions of Proposition \ref{6equicond} should then be
evaluated on
$\bar{L}_\mu$ equipped with the Cartan map $\bar{\mu}$\,.

\smallskip
As already mentioned above, a canonical
interpretation of an arbitrary complete lattice as a property lattice
can be realized by
taking a copy of $L\setminus\{0\}$ as states with the Cartan map defined by
$\mu(a):=]0,a]$, but in general this solution violates
superpositional faithfullness\,.
Clearly, besides atomistic lattices there are many examples that
do allow a superpositionally faithful interpretation,
in particular when generalizing to weak Cartan maps. But can we provide
such an interpretation
for any complete lattice, or, for any bounded meet-semilattice?
So we leave the following questions open\,:
\bqu
Does any complete lattice $L$ admits an interpretation as a property
lattice where the distributive
hull can be interpreted as the disjunctive hull, i.e., does there
exists a weak Cartan map with
as extension a Cartan map that is superpositionally faithful such
that the restriction of
${\cal D}_{\bar{\mu}}(\bar{L})$ to $L$-disjunctions is a distributive
hull of $L$?
\equ
We end by investigating which properties with respect to existing
suprema are preserved in this
passage from finite to infinitary conjunctions.
\bpr
For $L$ a bounded meet-semilattice and $\mu:L\to{\cal P}(\Sigma)$ a
weak Cartan map we have the
following\,:

\smallskip\noindent
(i) $\bigvee_L A=\bigvee_{\bar{L}_\mu}\!A$ whenever $\bigvee_L A$ exists\,;

\smallskip\noindent
(ii) If $\bigvee_L A$ is disjunctive then it is distributive\,;

\smallskip\noindent
(iii) If $\bigvee_L A$ is disjunctive then $\bigvee_{\bar{L}_\mu}\!A$
is disjunctive\,;

\smallskip\noindent
(iv) If $\bigvee_L A$ is disjunctive then $\bigvee_{\bar{L}_\mu}\!A$
is distributive.
\epr
\bpf
(i): Clearly, $\bigvee_L A\geq\bigvee_{\bar{L}_\mu}\!A$\,. We
moreover have that
\beqa
\bigvee_{\bar{L}_\mu}A\!\!
&=&\!\!\bigwedge_{\bar{L}_\mu}\{b\in\bar{L}_\mu|\forall a\in A:a\geq b\}\\
&\stackrel{\cong}{\leftrightarrow}&\!\!\bigcap\{B\in{\cal
C}_\mu(L)|\forall C\in\mu[A]:B\supseteq C\}\\
&\stackrel{\cong}{\leftrightarrow}&\!\!\bigcap\Bigl\{B\in\bigl\{\bigcap\mu[C]\bigm|C\subseteq
L\bigr\}
\Bigm|B\supseteq\bigcup\mu[A]\Bigr\}\\
&\stackrel{\cong}{\leftrightarrow}&\!\!\bigcap\Bigl\{B\in\mu[L]
\Bigm|B\supseteq\bigcup\mu[A]\Bigr\}\\
&=&\!\!\bigwedge_{\bar{L}_\mu}\{b\in L|\forall a\in A:b\geq a\}
\eeqa
and that for $b\in L$, if $\forall a\in A:b\geq a$ then
$b\geq\bigvee_L A$, so we also have
$\bigvee_{\bar{L}_\mu}\!A\geq\bigvee_LA$\,. (ii): Straightforward
verification along the lines of
Proposition
\ref{propDisThenDistr}\,. (iii):
$\bar{\mu}(\bigvee_{\bar{L}_\mu}A)=\mu(\bigvee_{\bar{L}_\mu}A)=\mu(\bigvee_LA)=
\bigcup\mu[A]=\bigcup\bar{\mu}[A]$\,.
(iv): Follows from (iii) and Proposition \ref{propDisThenDistr}\,.
\epf
\bqu
Do their exist and what are the explicit analogues of Proposition
\ref{6equicond} and Proposition \ref{DILas BOUND}
for weak Cartan maps?
What are the answers to the analogues of Question \ref{quortho} and Question
\ref{questionM(L)} for weak Cartan maps?
\equ
In this section we have briefly discussed a situation where we might
have infinitary non-conjunctive
infima, nor did we assume atomisticity.  Their are however other
situations considered in physics
where with an underlying complete and even atomistic property
lattice, e.g., classical or quantum
physics, one chooses to consider an incomplete non-atomistic subset,
for example for reasons imposed by the
very nature of measure theory which forces to restrict to
$\sigma$-completeness, e.g., Pt\'ak and
Pulmannov\'a (1991).
\bqu
To which extend do the constructions, representations,
interpretations and results of this paper hold, or
how should they be modified, when replacing completeness by
$\sigma$-completeness or when adopting the
settings of any approach within the general field of ordered quantum
structures.
\equ

\bigskip\noindent
{\bf 5. SUMMARY, CONCLUSION AND PERSPECTIVES}

\medskip\noindent
Since any complete lattice can canonically be embedded in a complete
Heyting algebra, where this
embedding itself equips the complete Heyting algebra with an
additional operation, and since that
whenever this complete lattice is the lattice of properties of a physical
system this complete Heyting algebra encodes the logical propositions
on these properties, we
are tempted to claim that quantum logic should not be seen as
contradicting intuitionism, but
{\it entailing a refinement of intuitionism encoded in terms of
operational resolution and operational
complementation}. Complete Heyting algebras saturate this embedding
into an isomorphism, encoding exactly
those property lattices where all logical expressions involving
disjunctions define themself a property of
the system, recalling here that suprema in property lattices are in
general not disjunctive but {\it
introduce superpositions} whereas infima are indeed conjunctive.
Since the Bruns-Lakser construction for
injective hulls in the category of meet-semilattices turns out to be
a distributive hull, it provides a
disjunctive hull for superpositionally faithful property lattices (either with respect to an
ordinary or a weak Cartan map).

\smallskip
It is our
feeling that the need to define actuality sets in order to encode
emergence of disjunctions in temporal
processes, e.g., measurements, is in a one to one way connected with
propositions on the system's
dynamical behavior. This claim is strengthened by the fact that
operational resolutions prove to be the
mathematical objects that naturally go equipped with state and
property transitions as morphisms (Coecke
and Stubbe 1999a,b). The fact that the considerations made in this
paper haven't been made before could
be connected to intrinsic static nature of what has been
conceived as quantum logic. However,
since these dynamical considerations formally encode in terms of
categories rather than in terms
of lattices and require a complementary conceptual discussion than
the one in this paper, we have chosen
to discuss the dynamical applications of the results of this paper in
a separate paper (Coecke
2001). 

\smallskip
Finally, the carefull (and probabely also the non-carefull) reader has noticed that nowhere
in the paper {\sl weak modularity} plays any role.  However, in Coecke and Smets (2001) the
claim is made that the transition from either classical or 
constructive/intuitionistic logic to quantum logic entails
besides the introduction of an
additional unary connective operational resolution the shift from a binary
connective implication to a ternary connective where two of the arguments have an
ontological connotation and the third, the new one, an empirical.
These ternary connectives have a fundamentally dynamic nature and have the intuitionistic
ones introduced in this paper as statical limit.  This second
aspect of the shift from classical or constructive/intuitionistic to quantum will then be
the one that requires {\sl orthomodularity} of the underlying lattice of properties as a
crucial feature.  

\bigskip\noindent
{\bf APPENDIX: IMPLICATION VIA ADJUNCTION}

\medskip\noindent
It is the aim of this paragraph to illustrate how one proceeded in
previous attempts to equip
quantum logic with an implication (Hardegree 1979, Kalmbach 1983). To the
present author's opinion, this can be expressed the best in terms of
adjointness between
action of conjunction and left action of the implication arrow.
Recall that\,:

\smallskip\noindent
(i) A pair of maps $f:L\to M$ and
$g:M\to L$ between posets $L$ and $M$ are {\sl Galois adjoint},
denoted by $f\dashv g$, if
and only if
$f(a)\leq b\Leftrightarrow a\leq g(b)$\,.

\smallskip\noindent
(ii) Whenever
$f\dashv g$, $f$ preserves existing suprema and $g$ existing infima.

\smallskip\noindent
(iii) For $L$ and
$M$ complete lattices, any {\it inf}-morphism $g:M\to L$
has a unique {\it sup}-preserving {\sl left adjoint}
$g_*:a\mapsto\bigwedge\{b\in M|a\leq
g(b)\}$ and any {\it sup}-morphism $f:L\to M$ a unique {\it
inf}-preserving {\sl right adjoint}
$f^*:b\mapsto\bigvee\{a\in L|f(a)\leq b\}$\,.

\smallskip\noindent
Setting $i:{\cal C}(L)\hookrightarrow L$ and
$i_*:L\to{\cal C}(L):a\mapsto{\cal C}(a)$ given a closure ${\cal C}$
on $L$\,, we have for $a\in L$ and
$b\in {\cal C}(L)$ that $a\leq b\Rightarrow{\cal C}(a)\leq{\cal
C}(b)=b$ and thus $i_*(a)\leq
b\Leftrightarrow {\cal C}(a)\leq b\Leftrightarrow a\leq
b\Leftrightarrow a\leq i(b)$ so
$i_*\dashv i$ where ${\cal C}=i\circ i_*$\,, i.e., any closure
factors in a {\it
sup}-endomorphism
$i_*$ and an {\it inf}-subobject inclusion $i$\,. Thus, the range
${\cal C}(L)$ of a closure
${\cal C}$ on a complete lattice $L$ is a {\it inf}-subobject of $L$,
and any {\it
inf}-subobject $M$ of $L$ defines a closure
${\cal C}_M:{L\to L}:a\mapsto{\bigwedge\{b\in M|a\leq b\}}$\,.
Notice that we have $(a\wedge-)\dashv(a\Rightarrow-)$ in any Heyting
semialgebra,
so $a\wedge-$ preserves existing joins what exactly results in saying
that the joins of all subsets
are distributive. Conversely,  if the supremum of every subset of a
meet-semilattice $H$
exists and is distributive, then $H$ is complete by Birkhoff's theorem, and
for all $a\in H$ the map
$a\wedge-:H\to H$ preserves all suprema so it has a unique right
adjoint $a\Rightarrow-:H\to
H$, as such encoding $(-\Rightarrow-)$ when viewing $a$ as an argument.
It then follows that complete Heyting algebras are complete lattices
where the suprema of all
subsets are distributive.
Now, recalling that a complete ortholattice $L$ is a {\sl complete
orthomodular lattice} if
it is moreover weakly modular, i.e., if $a\leq b$ implies
$a\vee(a'\wedge b)=b$\,, setting
$\varphi_a:L\to L:b\mapsto a\wedge(a'\vee b)$ and $\varphi_a^*:L\to L:b\mapsto
a'\vee(a\wedge b)$ we have
$\varphi_a\dashv\varphi_a^*$.  Indeed, if $a\wedge(a'\vee b)\leq c$
then $a'\vee\bigl(a\wedge
\left(a\wedge(a'\vee b)\bigr)\right)\leq a'\vee(a\wedge c)$ where
$b\leq a'\vee b=a'\vee\bigl(a\wedge(a'\vee b)\bigr)$ since $a'\leq
a'\vee b$\,, and analogously
one proves the converse. This adjunction embodies why $\varphi_{(-)}^*(-)$
has been interpreted as an implication, since
$\varphi_a$ coincides with $(a\wedge -):L\to L$ in the case that $L$
is distributive.
This view is moreover motivated by the fact that where for a Heyting
semialgebra the actions
$\{(a\wedge-)|a\in L\}$ can be envisioned as projections on $a$\,,
for orthomodular lattices  the {\sl
Sasaki projections} $\{\varphi_a|a\in L\}$ are the closed orthogonal
projections in the {\sl Baer
$^*$-semigroup  of $L$-hemimorphisms}\, (Foulis 1960)\,. For the
particular case of the lattice of closed
subspaces of a Hilbert space the action of these Sasaki
projections coincides with that of the
projection operators on the corresponding closed subspaces.
For details and a more general discussion on the matter we
respectively refer to Kalmbach (1983) and
Coecke, Moore and Smets (2001c)\,.

\bigskip\noindent
{\bf ACKNOWLEDGMENTS}

\medskip\noindent
I definitely should thank both the referee and John Harding for
pointing out the
existence of the Bruns and Lakser paper, restoring my (inadmissible)
ignorance concerning  --- a
first version of this paper which I made available as downloadable
postscript presented their
results as mine.  I moreover thank the referee taking so much care in
refereeing, resulting in a definite improvement of the readability and motivational content
of this paper, and, who actually
contributed this paper. If he wasn't anonymous, I should have
proposed him as co-author in
view of his twelve page constructive report.  We also thank John Harding, David
Moore, Sonja Smets, Isar Stubbe and Frank Valckenborgh for
additional comments.

\bigskip\noindent
{\bf REFERENCES}\,\footnote{Preprints and postscript files of
published papers by the
current author can be downloaded at http://www.vub.ac.be/CLEA/Bob/Coecke.html.}

\medskip\noindent
A{\scriptsize ERTS}, D. (1982) `Description of Many Separated
Physical Entities without the Paradoxes Encountered in Quantum
Mechanics', {\it Foundations of Physics} {\bf 12}, 1131.
\par
\vspace{2mm}
\par
\noindent
B{\scriptsize ANASCHEWSKI}, B. and B{\scriptsize RUNS}, G. (1967)
`Categorical
Characterization of MacNeille Completion', {\it Archiv der Mathematik}
{\bf 18}, 369.
\par
\vspace{2mm}
\par
\noindent
B{\scriptsize IRKHOFF}, G. and {\scriptsize VON} N{\scriptsize
EUMANN}, J. (1936) `The Logic of
Quantum Mechanics', {\it Annals of Mathematics} {\bf 37}, 823.
\par
\vspace{2mm}
\par
\noindent
B{\scriptsize RUNS}, G. and H{\scriptsize ARDING}, J. (2000)
`Algebraic Aspects of Orthomodular
Lattices', In\,: B. Coecke, D.J. Moore and A. Wilce, (Eds.), {\it
Current Research in Operational
Quantum Logic: Algebras, Categories and Languages},
pp.37--66, Kluwer Academic Publishers.
\par
\vspace{2mm}
\par
\noindent
B{\scriptsize RUNS}, G. and L{\scriptsize AKSER}, H. (1970)
`Injective Hulls of Semilattices',
{\it Canadian Mathematical Bulletin} {\bf 13}, 115.
\par
\vspace{2mm}
\par
\noindent
C{\scriptsize OECKE}, B. (2000) `Structural Characterization of
Compoundness', {\it
International Journal of Theoretical Physics} {\bf 39}, 581\,;
arXiv:quant-ph/0008054.
\par
\vspace{2mm}
\par
\noindent
C{\scriptsize OECKE}, B. (2001) `Disjunctive Quantum Logic in
Dynamic Perspective',
{\it Studia Logica} {\bf 71}, 1; arXiv: math.LO/0011209\,.  
\par
\vspace{2mm}
\par
\noindent
C{\scriptsize OECKE}, B., M{\scriptsize OORE}, D.J. and S{\scriptsize
METS}, S. (2001a) `From
Operationality to Logicallity I. Philosophical and Formal
Preliminaries',  Submitted.
\par
\vspace{2mm}
\par
\noindent
C{\scriptsize OECKE}, B., M{\scriptsize OORE}, D.J. and S{\scriptsize
METS}, S. (2001b) `From
Operationality to Logicallity II. Syntax and Semantics',  Submitted.
\par
\vspace{2mm}
\par
\noindent
C{\scriptsize OECKE}, B., M{\scriptsize OORE}, D.J. and S{\scriptsize
METS}, S. (2001c) `Adjoint
Implications in Lattice Logics',  Preprint.  
\par
\vspace{2mm}
\par
\noindent
C{\scriptsize OECKE}, B.  and S{\scriptsize METS}, S. (2001) `The Sasaki-Hook is not
a [Static] Implicative Connective but Induces a Backward [in Time] Dynamic One that
Assigns Causes', Paper submitted to {\it International Journal of
Theoretical Physics\,} for the proceedings of IQSA V, Cesena, Italy, April 2001;
arXiv:quant-ph/0111076\,.
\par
\vspace{2mm}
\par
\noindent
C{\scriptsize OECKE}, B. and S{\scriptsize TUBBE}, I. (1999a) `On a
Duality of Quantales
Emerging from an Operational Resolution', {\it International Journal
of Theoretical
Physics} {\bf 38}, 3269.
\par
\vspace{2mm}
\par
\noindent
C{\scriptsize OECKE}, B. and S{\scriptsize TUBBE}, I. (1999b)
`Operational Resolutions and
State Transitions in a Categorical Setting', {\it Foundations of
Physics Letters} {\bf
12}, 29\,; arXiv: quant-ph/0008020.
\par
\vspace{2mm}
\par
\noindent
F{\scriptsize OULIS}, D.J. (1960) `Baer $^*$-Semigroups', {\it
Proceedings of the American Mathematical
Society} {\bf 11}, 648.
\par
\vspace{2mm}
\par
\noindent
F{\scriptsize OULIS}, D.J. and R{\scriptsize ANDALL}, C.H. (1984) `A
Note on Misunderstandings
of Piron's Axioms for Quantum Mechanics', {\it Foundations of
Physics} {\bf 14}, 65.
\par
\vspace{2mm}
\par
\noindent
G{\scriptsize OLDBLATT}, R. (1984) `Orthomodularity is Not Elementary',
{\it Journal of Symbolic Logic} {\bf 49}, 401.
\par
\vspace{2mm}
\par
\noindent
H{\scriptsize ARDEGREE}, G.M. (1979) `The Conditional in Abstract and Concrete
Quantum Logic', In\,: C. Hooker, (Ed.),
{\it Logico-Algebraic Approach to Quantum Mechanics II}, pp.49--108,
Reidel Publishing
Company.
\par
\vspace{2mm}  
\par
\noindent
H{\scriptsize ORN}, A. and K{\scriptsize IMURA}, N. (1971) `The Category of
Semilattices', {\it Algebra Universalis} {\bf 1}, 26.
\par
\vspace{2mm}
\par
\noindent
J{\scriptsize AUCH}, J.M. and P{\scriptsize IRON}, C. (1969) `On the
Structure of Quantal
Proposition Systems', {\it Helvetica Physica Acta} {\bf 42}, 842.
\par
\vspace{2mm}
\par
\noindent
J{\scriptsize OHNSTONE}, P.T. (1982) {\it Stone Spaces}, Cambridge
University Press.
\par
\vspace{2mm}
\par
\noindent
K{\scriptsize ALMBACH}, G. (1983) {\it Orthomodular Lattices}, Academic Press.
\par
\vspace{2mm}
\par
\noindent
M{\scriptsize ALINOWSKI}, J. (1990) `The Deduction Theorem for
Quantum Logic --- Some Negative
Results', {\it Journal of Symbolic Logic} {\bf 55}, 615.
\par
\vspace{2mm}
\par
\noindent
M{\scriptsize OORE}, D.J. (1993) `Quantum Logic Requires Weak Modularity',
{\it Helvetica Physica Acta}
{\bf 66}, 471.
\par
\vspace{2mm}
\par
\noindent
M{\scriptsize OORE}, D.J. (1999)
`On State Spaces and Property Lattices',
{\it Studies in History and Philosophy of Modern Physics} {\bf 30}, 61.
\par
\vspace{2mm}
\par
\noindent
P{\scriptsize ASEKA}, J. (1994) `Covers in Generalized Frames', In\,:
Chajda, I.
(Ed.), {\it Proceedings of the International Conference and Summer School
on General Algebra and Ordered Sets 1994}, pp.84--99, Palack\'y
University Publishing, Olomouc.
\par
\vspace{2mm}
\par
\noindent
P{\scriptsize IRON}, C. (1976) {\it Foundations of Quantum Physics},
W.A. Benjamin, Inc.
\par
\vspace{2mm}
\par
\noindent
P{\scriptsize T\'AK}, P. and P{\scriptsize ULMANNOV\'A}, S. (1991)
{\it Orthomodular Structures as Quantum
Logics}, Kluwer Academic Publishers.
\par
\vspace{2mm}
\par
\noindent
S{\scriptsize TUBBE}, I. (2000)
`A Categorical View on Frame Completions of Meet-Semilattices by Means of Distributive
Joins', Privately communicated research notes.
\par
\vspace{2mm}  
\par
\noindent
V{\scriptsize ICKERS}, S. (1989) {\it Topology Via Logic}, Cambridge
University Press.

\end{document}